\documentclass[11pt]{amsart}
\usepackage{amsbsy,amssymb,amscd,amsfonts,latexsym,amstext,delarray,
amsmath,graphicx} 

\newtheorem{thm}{Theorem}[section]
\newtheorem{prop}[thm]{Proposition}
\newtheorem{cor}[thm]{Corollary}
\newtheorem{lem}[thm]{Lemma}
\newtheorem{defn}[thm]{Definition}

\numberwithin{equation}{section}

\def\bS{{\mathbb S}}

\def\C{{\mathbb C}}

\renewcommand{\H}{{\mathbb H}}
\def\N{{\mathbb N}}
\renewcommand{\P}{{\mathbb P}}
\def\Q{{\mathbb Q}}
\def\Z{{\mathbb Z}}
\def\R{{\mathbb R}}

\def\fA{{\mathfrak A}}
\def\fP{{\mathfrak P}}
\def\fS{{\mathfrak S}}

\def\cB{{\mathcal B}}

\def\cH{{\mathcal H}}

\def\cL{{\mathcal L}}

\def\cP{{\mathcal P}}

\def\cR{{\mathcal R}}

\def\cU{{\mathcal U}}
\def\cV{{\mathcal V}}
\def\cW{{\mathcal W}}

\title{Cuntz--Krieger algebras and wavelets on fractals}
\author{Matilde Marcolli}
\author{Anna Maria Paolucci}
\address{Department of Mathematics  \\
California Institute of Technology \\ 
Pasadena, CA 91125, USA}
\email{matilde\@@caltech.edu}
\address{Max Planck Institute for Mathematics\\
Vivatsgasse 7 \\ D-53111 Bonn, Germany}
\email{paolucci\@@mpim-bonn.mpg.de}

\begin{document}
\maketitle

\begin{abstract}We consider representations of Cuntz--Krieger algebras
on the Hilbert space of square integrable functions on the limit set,
identified with a Cantor set in the unit interval. We use these representations
and the associated Perron-Frobenius and Ruelle operators to construct 
families of wavelets on these Cantor sets.\end{abstract}

\tableofcontents

\section{Introduction}\label{Intro}

A class of representations of the Cuntz algebra $O_N$ called
permutative representations were studied and classified in
\cite{Da-Pi}, \cite{Br-Jo}, \cite{Br-Jo-Os}. Besides 
interest in their own right within the field of operator algebras,
Cuntz algebras representations have very interesting applications 
to wavelets, fractals, and dynamical systems, see \cite{Br-Jo}
and \cite{Br-Jo-Os}. Some of these results have been extended 
to the more general class of Cuntz-Krieger algebras (see \cite{Ka}, \cite{Ka2}, \cite{KeStaStr}),
where representations of these algebras are related to 
Perron--Frobenius operators of certain measure space transformations.
Similar representations of Cuntz--Krieger algebras were considered
in the context of limit sets of Schottky groups and actions on trees
in \cite{ConsMa1}, \cite{ConsMa2}, \cite{CMRV} for arithmetic applications
to Arakelov geometry and p-adic Mumford curves. 
 
In this paper we look at representations of the Cuntz--Krieger
algebra having a underlying self-similarity structure. 
The concept of self-similarity has proved to be fundamental in
mathematics as well as in diverse applications, related to the
renormalization of structures on nested families of scales. 
In the theory of wavelets, the scales may be represented in
resolutions taking the form of nested systems of linear spaces, while
in $C^*$-algebra theory it gives rise to representations of algebras
on generators and relations such as those that define the Cuntz and 
Cuntz--Krieger algebras.

Cuntz--Krieger algebras arise naturally from semibranching function
systems on measure spaces, where the partial inverses $\sigma_i$ of the coding
map $\sigma$ are not defined everywhere. The resulting algebra is generated
by partial isometries $S_i$ associated to the maps in the semibranching function
system, and the relations between these generators involve a matrix $A$ with
entries equal to zero or one, which describes the decomposition of the domains 
of the $\sigma_i$ as a union of ranges of other $\sigma_j$ in the same family.

Conversely, a Cuntz--Krieger algebra $O_A$ defined by generators and relations
in terms of an $N\times N$-matrix $A$ as above determines a semibranching function
system on the limit set of infinite sequences in an alphabet on $N$ letters with the
admissibility condition that consecutive letters $ij$ can appear in a word if and
only if the corresponding entry in the matrix is $A_{ij}=1$. One can identify this limit
set $\Lambda_A$ as a Cantor set inside the interval $[0,1]$ by considering points whose $N$-adic
digital expansion satisfies the admissibility condition. On this Cantor set the action
of the maps $\sigma_i$ become simple shifts in the $N$-adic expansion and 
the representation of $O_A$ on the Hilbert space $L^2(\Lambda_A,\mu)$, with
respect to the Hausdorff measure of the appropriate dimension, has an
especially simple form, and so does also the Perron-Frobenius operator
for the shift map $\sigma$, which is expressed in terms of the generators of
the algebra. 

The Hausdorff dimension of the limit set $\Lambda_A$ is computed using
the Perron--Frobenius theorem for the non-negative matrix $A$, which
also shows that the components of the Perron--Frobenius eigenvector 
of $A$ give the measures of the ranges of the maps $\sigma_i$ in the
normalized Hausdorff measure of dimension the Hausdorff dimension 
of $\Lambda_A$, which is the unique probability measure satisfying the 
self-similarity condition for the fractal set $\Lambda_A$. 

The Perron--Frobenius  eigenvector of the matrix $A^t$ determines
a fixed point for the Perron--Frobenius operator for the shift map $\sigma$ on
the limit set $\Lambda_A$, which in turn gives a KMS state for an associated
time evolution on the algebra $O_A$ at inverse temperature equal to the
Hausdorff dimension of $\Lambda_A$.

One can construct as in \cite{Jor2} further measures on $\Lambda_A$,
using operator valued measures and square-integrable functions
of unit norm. As in the case of the Cuntz algebras, by analyzing 
the Fourier transforms of these measures, one sees that one can
approximate them with Dirac measures supported at truncations of
the $N$-adic expansions.

Besides the Cantor set $\Lambda_A \subset [0,1]$, there is another
fractal set that one can associate to the same matrix $A$, namely a
Sierpinski fractal $\bS_A$ inside the unit cube, given by points $(x,y)$
whose digits in the $N$-adic expansion satisfy the condition that
$A_{x_i y_i}=1$. The Hausdorff dimension of these sets is simply
computed in terms of the number of non-zero entries in $A$.
The shifts in the $N$-adic expansion determine a semibranching function
systems on $\bS_A$, where, unlike in the case of $\Lambda_A$, the
maps are everywhere defined, hence they give rise to an action of
a Cuntz algebra of rank depending on the number of non-zero entries in $A$.
There is a natural embedding of $\Lambda_A$ into $\bS_A$ induced
by the shift map on $\Lambda_A$. The action of the Cuntz algebra
determines via this embedding of $\Lambda_A$ into $\bS_A$ 
an action of a Cuntz--Krieger algebra.

We show how to use the representation of the algebra $O_A$
to construct an orthonormal system of wavelets on $L^2(\Lambda_A,d\mu)$.

We then consider the Ruelle transfer operator for the shift
$\sigma$ on $\Lambda_A$, with non-negative valued 
potential $W$ satisfying the Keane condition that the sum of
the values over preimages under $\sigma$ adds up to one. We show
that one can construct from these measures on $\Lambda_{A^t}$,
for the transpose matrix $A^t$, in terms of random walks where the
probabilities assigned to words of a given length in the alphabet
depends upon the values of the potential $W$. A simple example
of a potential satisfying the Keane condition is given in terms of
trigonometric functions.

The example of the continued fraction expansion on the
Hensley Cantor sets considered in \cite{Mar}, \cite{Mar2} is
described as an example where the general results of this
paper can be applied. 

As an application we also show how the technique we described
to construct wavelets on the Cantor sets $\Lambda_A$ can be
adapted to construct families of {\em graph wavelets}, using Cuntz--Krieger
algebras associated to finite graphs with no sinks. Graph wavelets
are considered a useful tool for spatial network traffic analysis \cite{Crov}.

\section{Representations of Cuntz--Krieger algebras}\label{CUN}

Let $A$ be an $N\times N$ matrix $A$ with entries in $\{ 0,1 \}$. For
consistency with the notation we adopt later in the paper, it is
convenient to index the entries $A=(A_{ij})$ with indices $i,j\in \{
0,\ldots,N-1\}$ instead of $\{1, \ldots, N\}$.

Recall that the Cuntz-Krieger algebra $O_A$ associated to such a
matrix $A$ is the $C^{\ast }$-algebra 
generated by $N$ (non zero) partial isometries $S_0,\ldots ,S_{N-1}$ 
satisfying the relations 
\begin{equation}\label{CKrel1}
S_i^\ast S_i=\sum_{j} A_{ij} S_jS_j^\ast 
\end{equation}
and
\begin{equation}\label{CKrel2}
\sum_{i=0}^{N-1} S_i S_i^\ast =1. 
\end{equation}
The algebra $O_A$ is uniquely determined by the relations 
\eqref{CKrel1} and \eqref{CKrel2} 
and it is linearly spanned by the set of $S_{\alpha}S_{\beta}^{\star}$ 
with words $\alpha$ and $\beta$ in $\{ 0,\ldots,N-1 \}$ 
with possibly different lenghts $\left|\alpha\right|$ and 
$\left|\beta\right|$, see \cite{CK}.

We are especially interested here in representations of $O_A$ 
as bounded operators on Hilbert spaces of the form $\cH=L^2( X ,\mu)$,
for $(X,\mu)$ a measure space. 
The representations we are interested in, which include the
cases of the arithmetic applications mentioned above, are all 
constructed in terms
of what we refer to as a {\em semibranching function system}, which
will be concretely realized in terms of a shift map on a Cantor-like 
fractal set and its partial inverses. 

\begin{defn}\label{semibranchdef}
Consider a measure space $(X,\mu)$ and a finite 
family $\{ \sigma_i \}_{i\in I}$, $\# I=N$,
of measurable maps $\sigma_i : D_i \to X$, defined
on measurable subsets $D_i \subset X$.
The family $\{\sigma_i\}$ is a semibranching 
function system if the following holds.
\begin{enumerate}
\item There exists a corresponding family $\{D_i\}_{i=1}^N$ 
of measurable subsets of $X$ with the property that
\begin{equation}\label{sigmaiDi}
\mu(X\smallsetminus \cup_i R_i) =0, \ \ \ \text{ and } \ \ \
\mu(R_i\cap R_j) =0, \ \ \text{ for } \ \  i\neq j,
\end{equation}
where we denote by $R_i$ the range $R_i=\sigma_i(D_i)$. 
\item There is a Radon--Nikodym derivative 
$$ \Phi_{\sigma_i} = \frac{d (\mu\circ\sigma_i)}{d\mu} $$
with $\Phi_{\sigma_i} >0$, $\mu$-almost everywhere on $D_i$.
\end{enumerate}
A measurable map $\sigma: X\to X$ is called a coding map for
the family $\{ \sigma_i \}$ if $\sigma\circ \sigma_i(x) = x$
for all $x\in D_i$.
\end{defn}

Thus, the maps of the semibranching function system are partial
inverses of the coding map $\sigma$. Notice that the reverse
composition $\sigma_i \circ \sigma$ is only defined when the
image of $x$ under $\sigma$ lands in the domain $D_i$ of $\sigma_i$. 

\smallskip

Given a semibranching function system $\{\sigma_i\}_{i=0}^{N-1}$ with coding 
map $\sigma$, one can construct an associated family of linear operators
$\{T_i\}_{i=0}^{N-1}$ acting on the Hilbert space $L^2(X,\mu)$ by setting
\begin{equation}\label{Tiops}
(T_i\psi)(x)=\chi_{R_i}(x)\left(\Phi_{\sigma_i}(\sigma(x))\right)^{-1/2} 
\psi(\sigma(x)),
\end{equation}
with $\psi\in L^2(X,\mu)$, where $\chi_{R_i}$ is the characteristic 
function of $R_i\subset X$. 

\begin{lem}\label{TistarLem}
The adjoint of the operator $T_i$ of \eqref{Tiops} is of the form
\begin{equation}\label{Tistar}
(T_i^\star \xi)(x) = \chi_{D_i}(x) (\Phi_{\sigma_i}(x))^{1/2}
\xi(\sigma_i(x)). 
\end{equation}
\end{lem}

\proof We have
$$ \langle T_i \psi, \xi \rangle = \int_{R_i} \overline{
\left(\Phi_{\sigma_i}(\sigma(x))\right)^{-1/2}  
\psi(\sigma(x))} \xi(x) \, d\mu(x) $$
$$ = \int_{D_i} \overline{ \left(\Phi_{\sigma_i}(u)\right)^{-1/2} \psi(u)}
\xi(\sigma_i(u)) \, \frac{d\mu\circ\sigma_i}{d\mu}\, d\mu(u) $$ $$ =
\int_{D_i} \overline{\psi(u)}
\left(\Phi_{\sigma_i}(u)\right)^{1/2} \xi(\sigma_i(u))
 \, d\mu(x) =\langle \psi, T_i^* \xi \rangle, $$
where we used the fact that the $\Phi_{\sigma_i}$ are positive 
real valued. This gives \eqref{Tistar}.
\endproof

We then see easily that the operators $T_i$ and $T_i^*$ satisfy
the following relation.

\begin{prop}\label{TiCKrel2lem}
The operators $T_i$ of \eqref{Tiops} and their adjoints \eqref{Tistar}
satisfy the relations $T_i T_i^* = P_i$, where $P_i$ is the projection
given by multiplication by $\chi_{R_i}$. This gives $\sum_i T_i T_i^*
=1$. Similarly, $T_i^* T_i=Q_i$, where $Q_i$ is the projection given
by multiplication by the characteristic function $\chi_{D_i}$.
\end{prop}

\proof We write explicitly the action of the operator $T_i T_i^*$ on
elements $\xi\in L^2(X,d\mu)$. We have
$$ (T_i T_i^* \xi)(x) = \chi_{R_i}(x) \chi_{D_i}(\sigma(x))
\Phi_{\sigma_i}^{-1/2} (\sigma(x)) \Phi_{\sigma_i}^{1/2}(\sigma(x))
\xi(x) = \chi_{R_i}(x) \xi(x). $$
Equivalently, we can write
$$ \langle T_i^* \xi, T_i^* \xi \rangle =
\int_{D_i} \Phi_{\sigma_i}(x) \, |\xi(\sigma_i(x))|^2 \, d\mu(x) $$
$$ = \int_{R_i} \Phi_{\sigma_i}(\sigma(u)) \, |\xi(u)|^2 \,
\frac{d\mu\circ \sigma}{d\mu} \, d\mu(u). $$
Notice then that one has
\begin{equation}\label{dmusigma}
 \frac{d\mu\circ\sigma}{d\mu} |_{R_i} =
(\Phi_{\sigma_i}\circ\sigma)^{-1}, 
\end{equation}
so that we obtain
$$  \langle T_i^* \xi, T_i^* \xi \rangle = \int_{R_i} |\xi(u)|^2 \,
d\mu(u) =\langle P_i \xi, P_i \xi \rangle, $$
which gives $T_i T_i^* = P_i$, the
range projection on $L^2(X,\mu)$ realized by the multiplication 
operator by the characteristic function of the set $R_i$. By the
assumptions \eqref{sigmaiDi} on the semibranching function 
system we know that the projections $P_i$ are orthogonal and
that $\sum_i P_i = 1$.

We then consider the product $T_i^* T_i$. 
We have
$$ \langle T_i\xi, T_i\xi \rangle = \int_{R_i}
\Phi_{\sigma_i}^{-1}(\sigma(x))
\, | \xi (\sigma(x)) |^2 \, d\mu(x) $$
If $x\in R_i$ then $\sigma(x)\in D_i$ since $\sigma\circ\sigma_i =id$
on $D_i$. Thus, we write the above as
$$ \int_{D_i} \Phi_{\sigma_i}^{-1}(u) \, |\xi (u)|^2 \, \left(
\frac{d\mu\circ\sigma}{d\mu} \right)^{-1} \, d\mu(u) =
\int_{D_i} |\xi (u)|^2 \, d\mu(u) =\langle \chi_{D_i} \xi,\chi_{D_i}
\xi \rangle, $$ 
where we used again \eqref{dmusigma}. This gives $T_i^* T_i=Q_i$,
where $Q_i$ is the domain projection given by multiplication by
the characteristic function $\chi_{D_i}$. Unlike the range projections
$P_i$, the domain projections $Q_i$ are, in general, not orthogonal.
\endproof

When the maps $\sigma_i$ are defined everywhere on $X$, one obtains
from the operators $T_i$ and $T_i^*$ a representation of the Cuntz
algebra $O_N$ in the following way.

\begin{prop}\label{ONrep}
Let $\{\sigma_i\}$ be a semibranching 
function system on $X$, where the $\sigma_i$ are defined on all of
$X$, that is, $D_i =X$ for all $i=0,\ldots,N-1$. Then the operators
$T_i$ define a representation of the Cuntz algebra $O_N$ on the
Hilbert space $\cH=L^2(X,\mu)$. Namely, they satisfy the relations
\begin{equation}\label{CuntzRel}
T_i^* T_i =1, \ \ \ \sum_i T_i T_i^* =1.
\end{equation}
\end{prop}

\proof Under the assumption that the semibranching function system 
has $D_i=X$ for all $i\in I$, we obtain from Proposition \ref{TiCKrel2lem} 
above that the operators $T_i$ and $T_i^*$ of \eqref{Tiops} 
and \eqref{Tistar} satisfy $T_i^* T_i =1$. Moreover, we know
from Proposition \ref{TiCKrel2lem} that $T_i T_i^* = P_i$, the
range projections given by multiplication
by the characteristic functions $\chi_{R_i}$. Since these range
projections are orthogonal and the union of the $R_i$ exhausts $X$ up
to sets of measure zero, we obtain that $\sum_i T_i T_i^* =1$.
\endproof

In the case where the maps $\sigma_i$ are not defined everywhere on
$X$, but only on smaller domains $D_i\subset X$, one 
can then use the operators $T_i$ and $T_i^*$ of \eqref{Tiops} and
\eqref{Tistar} to construct representations of Cuntz--Krieger
algebras, when the domains $D_i$ have the property that
\begin{equation}\label{DiRjA}
\chi_{D_i} = \sum_j A_{ij} \chi_{R_j}.
\end{equation} 
The examples considered in \cite{ConsMa1}, \cite{ConsMa2}, \cite{CMRV}
are particular cases of this general procedure.

\begin{prop}\label{DiRjAreps}
Let $\{\sigma_i\}$ be a semibranching 
function system on $X$, where the $\sigma_i$ are defined on
subsets $D_i\subset X$ satisfying \eqref{DiRjA} (possibly up to sets
of measure zero). Also assume that $A_{ii}=1$ for all $i=0,\ldots,N-1$.
Then the operators $T_i$ and $T_i^*$ of
\eqref{Tiops} and \eqref{Tistar} satisfy the Cuntz--Krieger relations
\eqref{CKrel1} and \eqref{CKrel2}, namely
\begin{equation}\label{TiCKrels}
\sum_i T_i T_i^* = 1 \ \ \ \text{ and } \ \ \  T_i^* T_i = \sum_j A_{ij} T_i T_i^*,
\end{equation}
hence they determine a representation of the Cuntz--Krieger algebra $O_A$
on the Hilbert space $\cH=L^2(X,\mu)$.
\end{prop}

\proof Using \eqref{Tistar} and \eqref{DiRjA} we have
$$ (T_i^* \xi)(x) = \sum_j A_{ij}\, \chi_{R_j}(x) \,
\Phi_{\sigma_i}^{1/2}(x) \, \xi(\sigma_i(x)). $$
We then obtain
$$ (T_i T_i^* \xi)(x) = \sum_j A_{ij} \chi_{R_i}(x) \chi_{R_j}(\sigma(x)) \,
\Phi_{\sigma_i}^{-1/2}(\sigma(x))\, \Phi_{\sigma_i}^{1/2}(\sigma(x))\,
\xi(x) $$
$$ = \sum_j A_{ij} \chi_{R_{ij}}(x) \xi(x) = P_i \xi (x), $$
since we have from \eqref{DiRjA} that
$$ \cup_{j:\, A_{ij}=1} R_{ij}=\{ x\in R_i\,|, \sigma(x)\in
D_i\}=R_i. $$
Since the projections $P_i$ are orthogonal, we then obtain 
$$ \sum_i T_i T_i^* = 1 . $$
This gives \eqref{CKrel2} with $S_i=T_i$. Similarly, we have
$$ T_i^* T_i = Q_i $$
from Proposition \ref{TiCKrel2lem}, where $Q_i$ is the projection
given by multiplication by $\chi_{D_i}$. Using again \eqref{DiRjA}
this then gives
$$ T_i^* T_i = \sum_j A_{ij} P_j = \sum_j T_i T_i^*, $$
which gives \eqref{CKrel1} with $S_i=T_i$.
\endproof

We describe below an important special case of semibranching 
function system, which gives rise to representations of Cuntz--Krieger
algebras of the type described in Proposition \ref{DiRjAreps}. 

\subsection{Cantor sets and subshifts of finite
type}\label{subshiftSec}

Let $\fA$ be an alphabet in $N$ letters, which we can identify with
the set $\{ 0, \ldots, N-1\}$. Let $\Lambda_A$ be the set of all
infinite {\em admissible} words in the alphabet $\fA$, where the
admissibility condition is specified by an $N\times N$ matrix $A$ with
entries in $\{ 0,1 \}$. Namely, 
\begin{equation}\label{LambdaA}
\Lambda_A:=\{ w=\{ x_n\}_{n=0,1,\ldots} \,|\, x_i \in \fA, \,\,
A_{x_i,x_{i+1}}=1 \}. 
\end{equation}
We assume further that the matrix $A$ has the property that $A_{ii}=1$
for each $i=0,\ldots,N-1$, that is, that arbitrarily long strings made
of the same letters are allowed in the words of $\Lambda_A$.

The set $\Lambda_A$ can be topologized as a Cantor set, for example by
identifying it with the subset of the interval $[0,1]$ of numbers
whose base $N$ expansion satisfies the admissibility condition. 
However, Notice that, when we choose to view $\Lambda_A$ as a subset
of the interval $[0,1]$, which is convenient in what follows, we 
identify the rational numbers infinite periodic sequences rather
than with a finite $N$-adic expansion, so as to be able to act with
the shift map $\sigma$ on all of $\Lambda_A$.

Let $\delta_A$ be the Hausdorff dimension of the set $\Lambda_A$,
realized as a subset of the interval $[0,1]$ in this way. We can then
consider the Hilbert space $L^2(\Lambda_A,\mu_A)$, where $\mu_A$ is
the Hausdorff measure in the dimension $\delta_A$. 

We consider on $\Lambda_A$ the self-map given by the one-sided shift 
\begin{equation}\label{shift}
\sigma: \Lambda_A \to \Lambda_A, \ \ \ \sigma(x_0x_1x_2\ldots
x_n\ldots)=x_1x_2\ldots x_n \ldots 
\end{equation}

\begin{prop}\label{sigmaAprop}
The shift $\sigma$ is the coding map of
the semibranching function system
\begin{equation}\label{sigmaiA}
\sigma_i: D_i \to R_i, \ \ \  \sigma_i(w)=iw,
\end{equation}
where 
\begin{equation}\label{DisigmaA}
D_i=\{ w=\{ x_k \}\in \Lambda_A \,|\, A_{i,x_0}=1 \} 
\end{equation}
and
\begin{equation}\label{RisigmaA}
R_i=\{ w=\{ x_k \}\in \Lambda_A \,|\, x_0 =i \}=:\Lambda_A(i) .
\end{equation}
\end{prop}

\proof We show that the maps of \eqref{sigmaiA} form a
semibranching function system. We have
$$ \Lambda_A = \cup_i R_i, \ \ \ \text { with } \ \ \ R_i\cap R_j
=\emptyset, \ \ i\neq j, $$
hence the condition \eqref{sigmaiDi} of a semibranching function
system is satisfied. Moreover, the Radon-Nikodym derivative 
$$ \Phi_{\sigma_i}=\frac{d \mu_A\circ \sigma_i}{d\mu_A} $$
is well defined and positive, since the map $\sigma_i$ on
$\Lambda_A \subset [0,1]$ is realized by contractions and
translations. In fact, we can write the domain $D_i$ of the
map $\sigma_i$ as
\begin{equation}\label{DiAijRj}
 D_i =\cup_{j:\, A_{ij}=1} R_j.
\end{equation}
On each $R_j$ the map $\sigma_i$ is the restriction of the map of 
the $I_j\subset [0,1]$,
$$ I_j =\{ w\in [0,1]\,|\, x_0=j \}, $$
where $x_0$ is the first digit in the $N$-adic expansion of
$w=0.x_0x_1x_2\ldots$, that maps it to the subset $I_{jj}$ of
elements with first and second digit equal to $j$ composed
with a translation that maps isometrically $I_{jj}\to I_{ij}$
to the interval of all numbers with first digit $i$ and second digit
$j$. It is then clear that the shift map \eqref{shift} is a coding
map for this semibranching function system, since on each $D_i$ 
we have $\sigma\circ \sigma_i(w)=w$.
\endproof

One then sees easily that this gives a representation of the
Cuntz--Krieger algebra $O_A$ of the type described in Proposition
\ref{DiRjAreps} above.

\begin{prop}\label{LambdaArep}
The operators $T_i$ and $T_i^*$ of \eqref{Tiops} and \eqref{Tistar}
acting on $\cH_A=L^2(\Lambda_A,\mu_A)$ define a representation of
$O_A$ with generators $S_i=T_i$. 
\end{prop}

\proof The result immediately follows from Proposition
\ref{DiRjAreps}, upon noticing that the condition \eqref{DiAijRj}
is the needed relation \eqref{DiRjA}. We are assuming $A_{ii}=1$
for all $i$, so the hypothesis of Proposition \ref{DiRjAreps} are
satisfied. 
\endproof

It is well known (see \cite{CK}) that the abelian $C^*$-algebra
$C(\Lambda_A)$ sits naturally inside the Cuntz--Krieger algebra 
$O_A$ as the $C^*$-subalgebra generated by the range projections
$$ S_{x_1}\cdots S_{x_n} S_{x_n}^* \cdots S_{x_1}^*, $$
for arbitrary $x_i\in \fA$ and arbitrary $n$.

\subsection{Perron--Frobenius operator}\label{PFsec}

Consider the operator $T_\sigma: L^2(X,\mu) \to L^2(X,\mu)$ that
composes with the coding map $\sigma: X\to X$,
\begin{equation}\label{Tsigma}
(T_\sigma \psi)(x) = \psi(\sigma(x)).
\end{equation}
It is well known in the theory of dynamical systems that one
can associate to a self  map $\sigma: X\to X$ of a measure
space its Perron--Frobenius operator $\cP_\sigma$. This is 
defined as the adjoint of the composition \eqref{Tsigma} by
\begin{equation}\label{PFdef}
\int  \overline{\psi} \, \cP_\sigma (\xi) d\mu = \int  \overline{T_\sigma(\psi)}\,\xi \, d\mu.
\end{equation}

\begin{prop}\label{PFsigma}
Let $\{\sigma_i\}_{i=1}^N$ be a semibranching function system
with coding map $\sigma: X \to X$. Then the Perron--Frobenius
operator $\cP_\sigma$ is of the form
\begin{equation}\label{PFsigmai}
(\cP_\sigma \xi)(x)= \sum_i \chi_{D_i}(x) \, \Phi_{\sigma_i}(x) \, 
\xi(\sigma_i(x)).
\end{equation}
\end{prop}

\proof
In the inner product of $\cH=L^2(X,\mu)$ we find
$$ \langle T\psi, \xi \rangle = \int_X \overline{\psi(\sigma(x))}\xi(x) \,
d\mu(x)  $$ $$ =
\sum_i \int_{D_i} \overline{\psi(u)} \xi(\sigma_i(u))
\frac{d(\mu\circ\sigma_i)}{d\mu}\, d\mu(u) = \langle \psi, \sum_i
\chi_{D_i}\Phi_{\sigma_i}\,\, \xi\circ \sigma_i  \rangle. $$
\endproof

Notice the similarity of the Perron--Frobenius operator $\cP_\sigma$
to the operators $T_i^*$ of \eqref{Tistar} above. In fact, using
\eqref{Tistar} and Proposition \ref{PFsigma}, we easily get the
following, which was observed already in \cite{Ka}.

\begin{cor}\label{PFTicor}
Let $\{\sigma_i\}_{i=1}^N$ be a semibranching function system
with coding map $\sigma: X \to X$. Then the Perron--Frobenius
operator $\cP_\sigma$ is of the form
\begin{equation}\label{PFTi}
\cP_\sigma = \sum_i \Phi_{\sigma_i}^{1/2}\, T_i^* .
\end{equation}
\end{cor}

Notice that, in some particular cases, the functions $\Phi_{\sigma_i}$
may be constant, in which case \eqref{PFTi} gives just
a linear combination of the operators $T_i^*$. For example, in the cases
considered in \cite{ConsMa2} and \cite{CMRV} the functions
$\Phi_{\sigma_i}$ are locally constant, while they are not in the case
considered in \cite{ConsMa1}.

In the case of representations as in Proposition \ref{DiRjAreps},
we can express the Perron--Frobenius operator in terms of
the partial isometries $S_i$ in the following way. 

\begin{prop}\label{PFSiprop}
Let $\{\sigma_i\}$ be a semibranching 
function system on $X$, where the $\sigma_i$ are defined on
subsets $D_i\subset X$ satisfying \eqref{DiRjA} (possibly up to sets
of measure zero). 
Then the Perron--Frobenius operator $\cP_\sigma$ is a function of the
adjoints $S_i^*$ of the generators of the Cuntz--Krieger algebra 
$O_A$ and the multiplication operators by the functions
$\Phi_{\sigma_i}^{1/2}$ by
\begin{equation}\label{PFSi}
\cP_\sigma = \sum_i \Phi_{\sigma_i}^{1/2}\, S_i^* .
\end{equation}
In the case where the $\Phi_{\sigma_i}^{1/2}$ are constant over $D_i$,
the operator $\cP_\sigma$ belongs to the algebra $O_A$.
\end{prop}

\proof The hypothesis are the same as in Proposition \ref{DiRjAreps},
hence we know that the generators $S_i$ of the Cuntz--Krieger algebra
$O_A$ in the representation on $L^2(X,\mu)$ are given by the $T_i$
of \eqref{Tiops}. Then \eqref{PFTi} gives \eqref{PFSi}.
The case where the $\Phi_{\sigma_i}^{1/2}$ are constant over $D_i$ then
follows immediately from \eqref{PFSi}, since $\cP_\sigma$ is then a
linear combination of the $S_i^*$.
\endproof

To avoid having to assume that the $\Phi_{\sigma_i}$ are constant in
the result above (although this will in fact be the case in the main
example we will be considering later), one can more conveniently work
with representations of the Cuntz--Krieger algebras on the Hilbert
space of half-densities, analogous to the representations of the Cuntz
algebra considered in \cite{Jor}.

Recall that the Hilbert space $\tilde\cH$ of half densities consists of elements
of the form $\psi (d\mu/d\lambda)^{1/2}$, where $\psi \in L^2(X,d\mu)$
and $\mu << \lambda$ with $d\mu/d\lambda$ the Radon--Nikodym
derivative, which $\lambda$-a.e.~positive. Elements are considered 
modulo $\lambda$-a.e. equivalence and the inner product is given by
\begin{equation}\label{densityprod}
\left\langle \psi \left(\frac{d\mu}{d\lambda}\right)^{1/2}, h
\left(\frac{d\nu}{d\lambda}\right)^{1/2} \right\rangle = \int_X \overline{f} \,
\left(\frac{d\mu}{d\lambda}\right)^{1/2} \, h \, \left(\frac{d\nu}{d\lambda}\right)^{1/2} \,
d\lambda .
\end{equation}
One often writes elements of $\tilde\cH$ with the notation $\psi
\sqrt{d\mu}$. 

Given a semibranching function system on $X$ satisfying \eqref{DiRjA},
we can construct representations of the Cuntz--Krieger algebra $O_A$
on the space of half-densities of $X$, in much the same way as we did
in Proposition \ref{LambdaArep} on the space $L^2(X,d\mu)$.

\begin{prop}\label{halfdenseRep}
Let $\{\sigma_i\}$ be a semibranching 
function system on $X$, where the $\sigma_i$ are defined on
subsets $D_i\subset X$ satisfying \eqref{DiRjA}, possibly up to sets
of measure zero. Let $\tilde\cH$ be the Hilbert space of
half-densities on $X$. Consider the operators
\begin{equation}\label{tildeSi}
\tilde S_i (\psi \sqrt{d\mu}) = \chi_{R_i}\,\, (\psi\circ \sigma) \,\,
\sqrt{d\mu \circ \sigma}.
\end{equation}
These operators define a representation of the Cuntz--Krieger algebra $O_A$. 
\end{prop}

\proof To compute the adjoints $S_i^*$ we check
$$ \langle \tilde S_i (\psi \sqrt{d\mu}) , \xi \sqrt{\nu}\rangle =
\int_{R_i} \overline{\psi(\sigma(x))} \xi(x) \,
\left(\frac{d\mu(\sigma(x))}{d\lambda}\right)^{1/2}
\left(\frac{d\nu(x)}{d\lambda}\right)^{1/2} \, d\lambda(x) $$
$$ = \int_{D_i} \overline{\psi(u)} \xi(\sigma_i(u))\, \left(\frac{d\mu(u)}{d\lambda\circ\sigma_i}\right)^{1/2}
\left(\frac{d\nu(\sigma_i(u))}{d\lambda\circ\sigma_i}\right)^{1/2}
\frac{d\lambda\circ\sigma_i(u)}{d\lambda} \, d\lambda(u) $$
$$ = \int_{D_i} \overline{\psi(u)}
\xi(\sigma_i(u))\,\left(\frac{d\mu(u)}{d\lambda}\right)^{1/2}
\left(\frac{d\nu(\sigma_i(u))}{d\lambda}\right)^{1/2} \, d\lambda(u)
$$
$$ = \langle \psi \sqrt{d\mu}, \chi_{D_i}\, \xi\circ\sigma_i\,
\sqrt{d\nu\circ\sigma_i}\rangle, $$
which gives
\begin{equation}\label{tildeSistar}
\tilde S_i^*(\xi \sqrt{d\nu}) = \chi_{D_i}\,\, (\xi\circ\sigma_i) \,\,
\sqrt{d\nu\circ\sigma_i}.
\end{equation}
We then check that the operators $\tilde S_i$ and $\tilde S_i^*$
satisfy the Cuntz--Krieger relations \eqref{CKrel1} and
\eqref{CKrel2}. We have
$$  \langle \tilde S_i^* (\xi \sqrt{d\nu}), \tilde S_i (\xi
\sqrt{d\nu})\rangle = \int_{D_i} |\xi(\sigma_i(x))|^2 \,
\frac{d\nu\circ\sigma_i}{d\lambda} \, d\lambda(x) $$
$$ = \int_{R_i} |\xi(u)|^2 \frac{d\nu}{d\lambda\circ\sigma}
\frac{d\lambda\circ\sigma}{d\lambda} \, d\lambda(u) =\int_{R_i}|\xi|^2\frac{d\nu}{d\lambda} \, d\lambda,  $$
which shows that $\tilde S_i  \tilde S_i^*=\tilde P_i$, the range
projection given by multiplication by the characteristic function of
$R_i$, so that the relation \eqref{CKrel2} is satisfied by the
orthogonality of the projections $\tilde P_i$
$$ \sum_i \tilde S_i  \tilde S_i^*= 1 . $$
We also have
$$ \langle \tilde S_i (\psi \sqrt{d\mu}), \tilde S_i (\psi
\sqrt{d\mu})\rangle = \int_{R_i} |\psi(\sigma(x))|^2
\frac{d\mu\circ\sigma}{d\lambda} \, d \lambda(x) $$
$$ =\int_{D_i} |\psi(u)|^2 \frac{d\mu}{d\lambda\circ\sigma_i}
\frac{d\lambda\circ\sigma_i}{d\lambda} \, d\lambda(u) = \int_{D_i}
|\psi(u)|^2 \frac{d\mu}{d\lambda} d\lambda(u), $$
which shows that $\tilde S_i^* \tilde S_i =\tilde Q_i$, where $\tilde
Q_i$ is the domain projection given by multiplication by the
characteristic function of $D_i$. Using the relation \eqref{DiRjA}
this then gives
$$ \tilde S_i^* \tilde S_i = \sum_j A_{ij} \tilde S_i \tilde S_i^* $$
which shows that \eqref{CKrel1} is satisfied.
\endproof

We then compute explicitly the Perron--Frobenius operator of the 
coding map $\sigma:X \to X$ acting on the space of half-densities.

\begin{prop}\label{PFhalfdense}
Let $\sigma: X \to X$ be the coding map of a semibranching 
function system as in Proposition \ref{halfdenseRep} above.
The Perron--Frobenius operator $\tilde\cP_\sigma$ on the Hilbert space
of half-densities is given by
\begin{equation}\label{tildePFsigma}
\tilde\cP_\sigma  = \sum_i \tilde S_i^*  ,
\end{equation}
where $\tilde S_i$ are the generators \eqref{tildeSi} of the
representation of the Cuntz--Krieger algebra $O_A$ on $\tilde\cH$.
\end{prop}

\proof The translation operator associated to the shift map $\sigma: X
\to X$ is acting on the space of half-densities by
$$ \tilde T_\sigma (\psi \sqrt{d\mu}) = \psi\circ \sigma \, \sqrt{d\mu
\circ \sigma}. $$
The Perron--Frobenius operator $\tilde\cP_\sigma$ on $\tilde\cH$ is
the adjoint 
$$ 
\langle \tilde T_\sigma (\psi \sqrt{d\mu}), \xi \sqrt{d\nu} \rangle
= \langle \psi \sqrt{d\mu}, \tilde\cP_\sigma (\xi \sqrt{d\nu}) \rangle.
$$
This gives
$$ \int_X \overline{\psi(\sigma(x))}
\left(\frac{d\mu\circ\sigma}{d\lambda} \right)^{1/2} \xi(x) \,
\left(\frac{d\nu}{d\lambda} \right)^{1/2} d\lambda(x) = $$
$$ \sum_i \int_{R_i} \overline{\psi(u)} \left(\frac{d\mu}{d\lambda\circ \sigma_i}
\right)^{1/2} \xi(\sigma_i(u))\,
\left(\frac{d\nu\circ\sigma_i}{d\lambda\circ\sigma_i} \right)^{1/2}
\frac{d\lambda\circ\sigma_i}{d\lambda} d\lambda(u) $$
$$ = \sum_i \int_{R_i} \overline{\psi(u)}
\left(\frac{d\mu}{d\lambda}\right)^{1/2} \, \xi(\sigma_i(u))
\left(\frac{d\nu\circ\sigma_i}{d\lambda}\right)^{1/2}  d\lambda(u), $$
which gives
$$ \tilde\cP_\sigma (\xi \sqrt{d\nu})=\sum_i \chi_{R_i}\,\,
(\xi\circ\sigma_i) \,\, \sqrt{d\nu\circ\sigma_i}, $$
which is \eqref{tildePFsigma}.
\endproof

For example, in the case of the Cuntz--Krieger algebras considered in
\cite{ConsMa1}, \cite{ConsMa2}, \cite{CMRV}, where the representation
comes from the action of a Schottky group $\Gamma$ on its limit set, the
generators $S_i$ are associated to a symmetric set of generators 
$\fA=\{ \gamma_1, \ldots, \gamma_g, \gamma_1^{-1}, \ldots, \gamma_g^{-1}
\}$ of a Schottky group of genus $g$, and the matrix $A$ of the
Cuntz--Krieger algebra has $A_{ij}=1$ for $|i-j|\neq g$ and zero
otherwise, corresponding to the admissibility of the infinite
sequences $w=a_0 a_1 a_2 \cdots$ of elements of $\fA$
parameterizing points in the limit set $\Lambda_\Gamma$, namely
that $a_{i+1}\neq a_i^{-1}$. In this particular class of examples, the
Perron--Frobenius operator of Proposition \ref{PFhalfdense} has the
form
$$ \tilde\cP_\sigma = \tilde S_{\gamma_1}^* + \tilde
S_{\gamma_1^{-1}}^* + \cdots + \tilde S_{\gamma_g}^* + \tilde
S_{\gamma_g^{-1}}^*. $$
This resembles closely a Harper operator for the group $\Gamma$, save
for the important difference that the operators associated to the
symmetric set of generators of $\Gamma$ here are partial isometries
and not unitaries as in the usual Harper operator. 

\subsection{Projection valued measures}\label{projmeasSec}

We recall how one constructs projection-valued measures using
subdivions of compact metric spaces and subdivisions of projections in
Hilbert spaces. (We follow the notation and terminology of \cite{Jor3}
for the standard material we recall.) We then show how this technique 
applies to the representations of Cuntz-Krieger algebras described above.

We begin by recalling the notion of partitions and $N$-adic systems of
partitions of a metric space.

\begin{defn}\label{AXsubdiv}
Let $(X, d)$ be a compact metric space. For subsets $A\subset X$, define the diameter as
\begin{equation}\label{diamA}
\left|A\right|:=\sup\{ d(x,y)\,|\, x,y \in A\}.
\end{equation}
A partition $\cP$ of $X$ is a family $\{A(i)\}_{i\in I}$, for a (finite) index set $I$, with the property that
\begin{enumerate}
\item $\bigcup_i A(i)=X$.
\item $A(i)\cap A(j)=\emptyset$, for $ i\neq j$.
\end{enumerate}
For a given $N\geq 2$, an $N$-adic system of partitions of $X$ is a
family (indexed by $k\in \N$) of partitions $\cP_k$ of
$X$ into Borel subsets $A_k(a)$, indexed by elements of $\fA^k$, where
$\fA=\{0, \ldots, N-1\}$ is the given alphabet on $N$ letters, with
the properties: 
\begin{enumerate}
\item $|A_k(a)| = O(N^{-ck})$, for some $c>0$.
\item Every $A_{k+1}(b)$, with $b\in \fA^{k+1}$, is contained in some
$A_k(a)$, for some $a\in \fA^k$.
\end{enumerate}
\end{defn}

We then recall the equally well known notion of partitions of
projections in Hilbert spaces.

\begin{defn}\label{ProjHilb}
Let $\cH$ be a complex separable Hilbert space. A partition of
projections in $\cH$ is a collection $\{ P(i) \}_{i\in I}$ of
projections $P(i)=P(i)^*=P(i)^2$ such that 
\begin{enumerate}
\item $P(i)P(j)=0$, for $i\neq j$.
\item $\sum_i P(i) =1$.
\end{enumerate}
An $N$-adic system of partitions of $\cH$ into projections is a family
of partitions into projections $\{ P_k(a) \}$ indexed by $a\in \fA^k$
such that, for every $P_{k+1}(a)$, there is some $b\in \fA^k$ with
$P_k(b) P_{k+1}(a) = P_{k+1}(a)$.  
\end{defn}

We also recall the notion of operator valued measure.

\begin{defn}\label{OpMeas}
Denote by $\cB(X)$ the collection of Borel subsets of a compact metric
space $X$. A positive operator-valued function $E:\cB(X)\to \cL(\cH)$
defined on $\cB(X)$ with values in bounded linear operators on a
Hilbert space $\cH$ is called a $\sigma$ additive measure if, given a sequence
$B_1,B_2\ldots,$ in $\cB(X)$, such that $B_i\cap B_j=\emptyset$ for
$i\neq j$, one has 
\begin{equation}\label{EBorel}
E\left(\bigcup_i B_i\right)=\sum_i E(B_i).
\end{equation}
An orthogonal projection valued measure is a positive operator-valued
measure as above satisfying:
\begin{enumerate}
\item $E(B)=E(B)^*=E(B)^2$, for all $B\in \cB(X)$.
\item $E(B_1)E(B_2)=0$ when $B_1\cap B_2 =\emptyset$.
\item $E(X)=1$, the identity on $\cH$.
\end{enumerate}
\end{defn}

Note that the values $E(B_i)$ in \eqref{EBorel} are positive
operators, so we take the summation on the right hand side of 
\eqref{EBorel} to be convergent in the strong operator topology.

We are interested here in a particular construction of $N$-adic
partitions, for the metric Cantor set $\Lambda_A$ defined in
\eqref{LambdaA} above. As above, we consider the alphabet
$\fA=\{ 0,\ldots, N-1\}$. For any $k\in \N$, we denote by
$\cW_{k,A}\subset \fA^k$ the finite set of all admissible words of
length $k$ in the alphabet $\fA$,
\begin{equation}\label{WkA}
\cW_{k,A}=\{ a=(a_1,\ldots,a_k) \in \fA^k\,|\, A_{a_i,a_{i+1}}=1,
\, i=1,\ldots,k \}. 
\end{equation}
We also denote by $\Lambda_{k,A}(a)$ the clopen subset of the Cantor
set $\Lambda_A$ given by all words that start with a given $a\in
\cW_{k,A}$, 
\begin{equation}\label{LambdakAa}
\Lambda_{k,A}(a)=\{ w=(w_1,w_2,\ldots,w_n,\ldots)\in \Lambda_A\,|\,
(w_1,\ldots,w_k)=a \}.
\end{equation}

We then have the following partition and corresponding operator valued
measure.

\begin{prop}\label{NadicCLambda}
The subsets $\Lambda_{k,A}(a)$ of \eqref{LambdakAa} define an $N$-adic
system of partitions for $\Lambda_A$. There is a correposnding
$N$-adic system of projections $P_k(a)$ on the Hilbert space
$\cH=L^2(\Lambda_A,\mu_A)$ and an orthogonal projection valued measure
$E$ on $\cB(\Lambda_A)$ satisfying
\begin{equation}
E(\Lambda_{k,A}(a)) = P_k(a),
\end{equation}
for all $k\in \N$ and for all $a\in \cW_{k,A}$.
\end{prop}

\proof To see that the $\Lambda_{k,A}(a)$ form an $N$-adic
system of partitions, notice that, when we identify $\Lambda_A$ with
the subset of $[0,1]$ of numbers with admissible $N$-adic digital
expansion and we measure diameters in the Euclidean distance on
$[0,1]$, we see that the set $\Lambda_{k,A}(a)$, which consists of
such numbers with fixed first $k$ digits in the $N$-adic expansion
have 
\begin{equation}\label{diamAbound}
 |\Lambda_{k,A}(a)| \leq N^{-k}, 
\end{equation}
since the sets of all numbers with fixed $k$ digits in the $N$-adic
expansion are intervals of length $N^{-k}$. Moreover, by construction
we have inclusions
\begin{equation}\label{inclkkmin1}
 \Lambda_{k,A}(a_1,\ldots,a_k)\subset
\Lambda_{k-1,A}(a_1,\ldots,a_{k-1}). 
\end{equation}
We also have, for fixed $k$, 
$$ \Lambda_{k,A}(a)\cap \Lambda_{k,A}(b)=\emptyset, \ \ \text{ for }
a\neq b \in \cW_{k,A}, $$
and 
$$ \cup_{a\in \cW_{k,A}} \Lambda_{k,A}(a) = \Lambda_A. $$
Thus, we have an $N$-adic system of partitions.

One knows from \cite{CK} that there is an $*$-isomorphism
between the $C^*$-algebra of continuous functions $C(\Lambda_A)$ 
and the maximal abelian subalgebra of the Cuntz--Krieger algebra 
$O_A$ generated by all the range projections
\begin{equation}\label{rangePka}
P_k(a)= S_{a_1} \cdots S_{a_k} S_{a_k}^* \cdots S_{a_1}^*.
\end{equation}
We show that the $P_k(a)$ define an $N$-adic system of projections 
on the Hilbert space $\cH=L^2(\Lambda_A,\mu_A)$. In the representation
of $O_A$ described in Proposition~\ref{DiRjAreps}, the operator
$P_k(a)$ acts as the projection given by multiplication by the
characteristic function of the set $\Lambda_{k,A}(a)$. 

Since the $\Lambda_{k,A}(a)$ form an $N$-adic system of partitions,
in particular, as we have seen above, there are inclusions
\eqref{inclkkmin1}. These imply that the corresponding projections
satisfy $$ P_{k-1}(a_1,\ldots,a_{k-1})P_k(a_1,\ldots,a_k) =
P_k(a_1,\ldots,a_k). $$

More precisely, one can see by writing as in \eqref{rangePka} and
using \eqref{CKrel2} that
$$ \sum_{a_k \in \fA} P_k(a_1,\ldots,a_k) = \sum_{a_k \in \fA} S_{a_1} \cdots S_{a_k} S_{a_k}^* \cdots
S_{a_1}^* $$
$$ = S_{a_1} \cdots S_{a_{k-1}} \left(\sum_{a_k \in \fA} S_{a_k} S_{a_k}^*\right) S_{a_{k-1}}^* \cdots
S_{a_1}^* $$ $$ = S_{a_1} \cdots S_{a_{k-1}} S_{a_{k-1}}^* \cdots
S_{a_1}^* =P_{k-1}(a_1,\ldots,a_{k-1}). $$

For every $k\in \bf Z_+$, let us denote by $\cU_k$ the finite
dimensional subalgebra of $C(\Lambda_A)$ spanned by the
finite linear combinations
$$ \sum_{a\in \cW_{k,A}} c_a \, \chi_{\Lambda_{k,A}(a)}. $$
The inclusions \eqref{inclkkmin1} determine embeddings
$\cU_{k-1} \to \cU_k$ and the bound \eqref{diamAbound} on the
diameters implies that every function in $C(\Lambda_A)$ can be
uniformly approximated with a sequence of functions in
$\cU=\varinjlim_k \cU_k$. Thus, the homomorphism
\begin{equation}\label{pihomo}
 \pi: \sum_{a\in \cW_{k,A}} c_a \, \chi_{\Lambda_{k,A}(a)} \mapsto
\sum_{a\in \cW_{k,A}} c_a \, P_k(a) 
\end{equation}
extends, by a standard argument from function theory, 
from $C(\Lambda_A)$ to all the Baire functions on $\Lambda_A$.

It makes sense then to define an operator valued measure by setting
\begin{equation}\label{ELambdaAmeas}
E(B):= \pi(\chi_B),
\end{equation}
where we still denote as $\pi$ the extension above. 
It follows that $E(\cdot)$ satisfies the properties of Definition
\ref{OpMeas} and is countably additive.
It also satisfies $E(\Lambda_{k,A}(a))=P_k(a)$, for every 
$k\in{\bf Z_+}$ and for all $a\in\cW_{k,A}$. 
\endproof

\subsection{Hausdorff dimension}\label{dimHsec}

We consider again the space $\Lambda_A$ of numbers in the interval
$[0,1]$ whose $N$-adic expansion is admissible according to the 
matrix $A$, that is, $x=0.a_0a_1\cdots a_n \cdots$ with
$A_{a_i,a_{i+1}}=1$. 

We know that in this case the maps $\sigma_i$ are defined on
domains $D_i \subset \Lambda_A$ satisfying $D_i = \cup_{j: A_{ij}=1}
R_j$, where $R_j\subset \Lambda_A$ is the range of $\sigma_j$, 
with $\Lambda_A=\cup_j R_j$ and $R_i\cap R_j=\emptyset$ when $i\neq
j$. We then have the following properties.

\begin{thm}\label{LambdaDim}
Assume that the non-negative matrix $A$ is irreducible, that is, there
exists a power $A^n$ for which all entries are positive.
Let $\delta_A$ be the Hausdorff dimension of $\Lambda_A$ and 
$\mu_A=\mu_{\cH,\delta_A}$ the corresponding Hausdorff measure.
\begin{enumerate}
\item On the sets $D_i\subset \Lambda_A$, the Radon-Nikodym derivatives 
are constant and equal to 
\begin{equation}\label{RNder}
\Phi_{\sigma_i}= \frac{d\mu\circ
\sigma_i}{d\mu} =N^{-\delta_A}.
\end{equation}
\item The Hausdorff measure $\mu=\mu_A$ on $\Lambda_A$ satisfies
\begin{equation}\label{muHRi}
\mu(R_i) = p_i,
\end{equation}
where $p=(p_i)_{i=0,\ldots,N-1}$ is the Perron--Frobenius eigenvector
of the matrix $A$, 
\begin{equation}\label{PFforA}
\sum_j A_{ij}\, p_j = r(A)\, p_i,
\end{equation}
with eigenvalue the spectral radius $r(A)$, and
normalized to have $\sum_i p_i =1$.
\item The Hausdorff dimension of $\Lambda_A$ is given by
\begin{equation}\label{dimHLambda}
\delta_A = \dim_H(\Lambda_A)= \frac{\log r(A)}{\log N},
\end{equation}
with $r(A)$ the spectral radius of the matrix $A$. 
\item The measure $\mu$ satisfies the self-similarity condition
\begin{equation}\label{musigmainv}
\mu = N^{-\delta_A} \sum_{k=0}^{n-1} \mu \circ \sigma_k^{-1},
\end{equation}
where $\mu(\sigma_k^{-1}(E))=\mu(\{ x\in \Lambda_A\,|\, \sigma_k(x)\in
E \})$.
\end{enumerate}
\end{thm}

\proof 
(1) The maps $\sigma_i$ act as the restrictions to the set $D_i$
of the linear maps
\begin{equation}\label{sigmaiLin}
\sigma_i (x) = \frac{x+i}{N}
\end{equation}
defined on the interval $[0,1]$. Thus, we see directly that the
Radon-Nikodym derivative of the Hausdorff measure $\mu_{\cH,s}$ 
will give
$$ \frac{d\mu_{\cH,s}\circ \sigma_i}{d\mu_{\cH,s}} = N^{-s}. $$
In particular for $s=\delta_A=\dim_H(\Lambda_A)$ this gives
\eqref{RNder}.

(2) We first show that setting 
\begin{equation}\label{nuHRi}
\nu(R_i) = p_i,
\end{equation}
with $p$ the normalized Perron--Frobenius eigenvector of $A$, 
defines a probability measure on $\Lambda_A$.

The Perron--Frobenius theorem for the matrix $A$ shows that, if
$r(A)$ denotes the spectral radius of $A$, then $r(A)$ is an
eigenvalue which has an eigenvector $p=(p_i)$ with non-negative
entries. We can normalize it so that $\sum_i p_i=1$. 
Setting $\nu(R_i)=p_i$ defines a measure on $\Lambda_A$. In fact,
it suffices to see that we can define $\nu(\Lambda_{k,A}(a))$
compatibly, for all $a\in \cW_{k,A}$. We set
\begin{equation}\label{nuLambdaka}
\nu(\Lambda_{k,A}(a))= r(A)^{-k} p_{a_k},
\end{equation}
where $a=(a_1,\ldots,a_k)\in \cW_{k,A}$. To see that
\eqref{nuLambdaka} consistently defines a measure on
$\Lambda_A$ we need to check that
\begin{equation}\label{nucompatible}
\nu(\Lambda_{k,A}(a)) = \sum_{j=0}^{N-1} A_{a_k j}\, \nu(\Lambda_{k+1,A}(aj)).
\end{equation}
We have
$$ \sum_j A_{a_k j}\, \nu(\Lambda_{k+1,A}(aj)) = \sum_j A_{a_k j}
r(A)^{-k-1} p_j = r(A)^{-k} p_{a_k} = \nu(\Lambda_{k,A}(a)), $$
where we used the Perron--Frobenius relation
$$ p_{a_k} = r(A)^{-1} \sum_j A_{a_k j} p_j. $$
The measure $\nu$ thus satisfies the self-similarity property
\begin{equation}\label{nusigmainv2}
\nu = r(A)^{-1} \sum_{j=0}^{n-1} \nu \circ \sigma_j^{-1}.
\end{equation}
Indeed, it suffices to check it on sets of the form
$E=\Lambda_{k,A}(a)$, for which $\sigma_j^{-1}(\Lambda_{k,A}(a))$ is
nonempty for $a_1=j$, in which case it is $\Lambda_{k-1,A}(\sigma(a))$.
Then we have
$$ \nu(\Lambda_{k,A}(jb))=r(A)^{-1} r(A)^{-|b|} p_{b_k} = r(A)^{-1}
\nu(\Lambda_{k-1,A}(\sigma(a))), $$
which gives \eqref{nusigmainv2}.

We then compare this with the Hausdorff measure $\mu=\mu_A$. This
satisfies
\begin{equation}\label{muNA}
 \mu(R_i) = N^{-\delta_A} \sum_j A_{ij} \mu(R_j). 
\end{equation}
In fact, this follows simply from the fact shown in (1) that
the Radon--Nikodym derivatives are constant,
$$ \Phi_{\sigma_i} = \frac{d\mu\circ \sigma_i}{d\mu} = N^{-\delta}, $$
which gives
$$ \mu(R_i)= \int_{D_i} \frac{d\mu\circ\sigma_i}{d\mu} d\mu =
N^{-\delta_A} \mu(D_i)= N^{-\delta_A}\sum_j A_{ij} \mu(R_j). $$
Note that it then follows that the measure $\mu$ also satisfies
\begin{equation}\label{muLambdaka}
\mu(\Lambda_{k,A}(a)) = N^{-k \delta_A} \mu(R_{a_k}),
\end{equation}
for $a=(a_1,\ldots,a_k)$. This follows directly from \eqref{muNA} and
the fact that
$$ \mu(\Lambda_{k,A}(a)) =\sum_j A_{a_k j}
\mu(\Lambda_{k+1,A}(aj)). $$
Notice then that \eqref{muNA} is saying that the vector $q=(q_i)$ with
$q_i =\mu(R_i)$ is also an eigenvector of the matrix $A$, with
eigenvalue $N^{-\delta_A}\leq r(A)$, with the normalization $\sum_i
q_i =1$. 

Under the assumption that the non-negative matrix $A$ is irreducible, the
Perron--Frobenius theorem for $A$ ensures that the eigenvalue $r(A)$
is simple and that if $q=(q_i)$ is another eigenvector, $A q=\lambda
q$ with $q_i \geq 0$, then $\lambda =r(A)$ and $q$ is a scalar
multiple of $p$. Since both vectors are normalized, this implies that 
\begin{equation}\label{PFsame}
N^{\delta_A}=r(A) \ \ \ \text{ and }  \ \ \ \nu(R_i)=p_i = q_i = \mu(R_i). 
\end{equation} 
By \eqref{muLambdaka} and \eqref{nuLambdaka}, this implies that
the measures $\mu$ and $\nu$ agree. 

(3) then follows immediately from $r(A)=N^{\delta_A}$ and (4) is just
the self-similarity \eqref{nusigmainv2}.
\endproof

As a particular case, if the matrix $A$ has the property that the
value $\alpha=\sum_j A_{ij}$ is the same for all $i=0,\ldots,N-1$, then one
has uniform probability for all the $R_i$, equal to $\mu(R_i)=1/N$, and
the set $\Lambda_A$ has then Hausdorff dimension $\delta_A=
\log(\alpha)/\log(N)$.

\medskip

We return to consider now in particular the representation of the
Cuntz--Krieger algebra $O_A$ on the space $L^2(\Lambda_A,d\mu_A)$ 
as in \S \ref{subshiftSec}. 

\begin{cor}\label{repLambdaAPF}
The Perron--Frobenius operator $\cP_\sigma$ on the Hilbert space
$L^2(\Lambda_A,d\mu_A)$, with $\mu_A=\mu_{\cH,\delta_A}$ the Hausdorff
measure with $\delta_A=\dim_H(\Lambda_A)$ satisfies
\begin{equation}\label{PFLambdaA}
\cP_\sigma = N^{-\delta_A/2} \sum_i S_i^*, 
\end{equation}
\end{cor}

\proof
As we have seen in Proposition
\ref{LambdaDim}, in this case the $\Phi_{\sigma_i}$ are locally
constant and equal to $N^{-\delta_A}$, with $\delta_A$ the Hausdorff
dimension, which in turn is given in terms of the spectral radius of
$A$. Then we have from Proposition \ref{PFSiprop} that the
Perron--Frobenius operator $\cP_\sigma$ on $L^2(\Lambda_A,d\mu_A)$ is
simply given by \eqref{PFLambdaA}, where the $S_i$ generate the 
representation of the Cuntz--Krieger algebra on $L^2(\Lambda_A,d\mu_A)$. 
\endproof

We then see that one can use the result of Theorem \ref{LambdaDim} to
construct a fixed point for the Perron--Frobenius operator
$\cP_\sigma$.

\begin{prop}\label{PFfix}
Assume that the matrix $A$ is irreducible, and let $\omega$ be the 
Perron--Frobenius eigenvector for $A^t$. Then $f=\sum_i \omega_i
\chi_{R_i}$ is a fixed point of the Perron--Frobenius operator
$\cP_\sigma$. 
\end{prop}

\proof Let $\omega$ be the Perron--Frobenius eigenvector
$$ A^t \omega = r(A) \omega . $$
The Perron--Frobenius operator $\cP_\sigma$ acting on the
function $f=\sum_i \omega_i \chi_{R_i}$ gives
$$ \cP_\sigma(f) = N^{-\delta_A} \sum_i \chi_{D_i} f\circ \sigma_i, $$
by Corollary \ref{repLambdaAPF}. We have
$$ \chi_{R_k}\circ \sigma_i =\delta_{ik} \chi_{D_i}, $$
which gives
$$ \cP_\sigma(f) = N^{-\delta_A} \sum_i \omega_i \chi_{D_i} =
N^{-\delta_A} \sum_{ij}  \omega_i A_{ij} \chi_{R_j} $$
from \eqref{DiRjA}. Using then $A^t \omega =r(A)\omega$
we obtain
$$ \cP_\sigma(f) = N^{-\delta_A} r(A) \sum_i \omega_i \chi_{R_i} =
\sum_i \omega_i \chi_{R_i} = f, $$
where we used the fact that $r(A)=N^{\delta_A}$ as in \eqref{PFsame}. . 
\endproof

\smallskip

There is a well known relation for Cuntz--Krieger
algebras between the fixed points of the dual Perron--Frobenius 
(or Ruelle tranfer operator) acting on measures and KMS states with 
respect to associated time evolutions, see \cite{KeStaStr}. We discuss
the more general case of the Ruelle transfer operators later, but we 
comment here on the case that follows directly from Theorem
\ref{LambdaDim}. 

\begin{cor}\label{KMSandPF}
On the Cuntz--Krieger algebra $O_A$ consider the time evolution
defined by setting 
\begin{equation}\label{NgaugeCK}
 \sigma_t(S_i) =N^{it} S_i. 
\end{equation}
The measure $\mu=\mu_A$ on $\Lambda_A$ 
defines a KMS state for the system $(O_A,\sigma_t)$
at inverse temperature $\beta = \delta_A$.
\end{cor}

\proof We define a state $\varphi$ on $O_A$ associated to the measure
$\mu$ by setting
\begin{equation}\label{stateSab}
 \varphi(S_a S_b^*) =\left\{ \begin{array}{ll} 0 & a\neq b \\[3mm] 
\mu(\Lambda_{k,A}(a)) & a=b \in \cW_{k,A}. 
\end{array}\right. 
\end{equation}
We use here the fact that all elements in $O_A$ can be approximated by
linear combinations of elements of the form $S_a S_b^*$. 
We then need to check that the state $\varphi$ satisfies the KMS
condition at inverse temperature $\beta=\delta_A$ for the time
evolution \eqref{NgaugeCK}. Because of the form of the state
\eqref{stateSab}, and the fact that the measure $\mu$ satisfies
\eqref{muLambdaka}, it suffices to check that
$$ \varphi(S_i^* S_i) = N^\beta \varphi(S_i S_i^*). $$
This follows since we have
$$ \varphi(S_i^* S_i)=\sum_j A_{ij} \varphi(S_j S_j^*)=
\sum_j A_{ij} \mu(R_j) = N^{\delta_A} \mu(R_i) = N^{\delta_A} \varphi(S_i S_i^*), $$
using the fact that $p=(p_i)$ with $p_i=\mu(R_i)$ is the
Perron--Frobenius eigenvector of the matrix $A$.
\endproof

\subsection{Real valued measures and Fourier transforms}

Given an element $f\in \cH$ with norm $\| f \| =1$, one can define a
real valued measure on $\Lambda_A\subset [0,1]$ by setting
\begin{equation}\label{mufE}
\mu_f(B) :=\langle f, E(B) f \rangle,
\end{equation}
with $E(B)$ an operator valued measure as in \S \ref{projmeasSec}.

Since each such $\mu_f$ is a compactly supported measure on the real
line, it makes sense to consider its Fourier transform
\begin{equation}\label{Fouriermu}
\widehat{\mu_f}(t) := \int e^{itx}\, d\mu_f(x) . 
\end{equation}

We then have the following result, which is analogous to
the case of the Cuntz algebras $O_n$ discussed in \cite{Jor2}.

\begin{prop}\label{FourierPF}
For every function $f\in \cH=L^2(\Lambda_A,d\mu_A)$ 
with $\| f \|=1$, the measure $\mu_f(E)=\langle f, P(E) f\rangle$ 
satisfies
\begin{equation}\label{mufSstar}
\sum_{k=0}^{N-1} \int_{\Lambda_A} \psi\circ \sigma_k \, d\mu_{S_k^* f} = \int_{\Lambda_A} \psi d\mu_f.
\end{equation}
The Fourier transform $\widehat \mu_f(t)$ satisfies
\begin{equation}\label{scaleFourier}
\widehat \mu_f(t) = \sum_{k=0}^{N-1}
e^{\frac{itk}{N}} \widehat{\mu_{S_k^* f}}(\frac{t}{N}).
\end{equation}
\end{prop}

\proof We have
\begin{equation}\label{mufSstar1}
\sum_k \int_{\Lambda_A} \psi\circ \sigma_k \, d\mu_{S_k^* f} =\sum_k \langle S_k^* f, 
\pi(\chi_{D_k}\,\psi \circ \sigma_k) S_k^* f \rangle, 
\end{equation}
where $\pi$ denotes the embedding $\pi: C(\Lambda_A) \hookrightarrow O_A$, 
as in \eqref{pihomo}, which realizes $C(\Lambda_A)$ as an abelian 
$*$-subalgebra of $O_A$, with $\pi(\chi_{\Lambda_{k,A}(a)})=S_a S_a^*$.
In the algebra $O_A$ we have the relations
\begin{equation}\label{CKfrel}
\begin{array}{l} 
\pi(f) S_k = S_k \, \pi(\chi_{D_k} \, f\circ \sigma_k) \\[2mm]
S_k \pi(f) = \pi(f\circ\sigma)\, S_k \\[2mm]
\pi(f) S_k^* = S_k^*\, \pi(f\circ\sigma) \\[2mm]
S_k^*\pi(f) =\pi(\chi_{D_k}\, f\circ \sigma_k)\, S_k^*.
\end{array} 
\end{equation}
Thus, we have $\pi(\chi_{D_k}\,\psi \circ \sigma_k) S_k^* = 
S_k^* \pi(\psi)$ and we write \eqref{mufSstar1} as
\begin{equation}\label{mufSstar12}
 \sum_k \langle f, S_k S_k^* \pi(\psi) f\rangle =
\sum_k \langle f, \pi(\chi_{R_k}\, \psi) f \rangle = \sum_k \int_{R_k}
\psi\, d\mu_f ,
\end{equation}
which gives \eqref{mufSstar}. We then proceed as in \cite{Jor2}, and
observe that \eqref{mufSstar}, applied to $\psi(x)=e^{itx}$, gives
$$ \sum_k \int e^{it\frac{x+k}{N}} \, d\mu_{S_k^* f}(x) =
\int e^{itx} \, d\mu_f(x) , $$
which gives \eqref{scaleFourier}.
\endproof

We can equivalently see \eqref{mufSstar} as an immediate consequence
of \eqref{musigmainv}, since we have
$$ \int \psi d\mu_f = \langle f, \pi(\psi) f\rangle =\int \psi |f|^2 d\mu $$
$$ =N^{-\delta} \sum_j \langle \chi_{D_j}\, f\circ\sigma_j ,
\pi(\psi\circ \sigma_j) \chi_{D_j}\, f\circ\sigma_j \rangle $$ $$ =
\sum_j \langle S_j^* f , \pi(\psi\circ \sigma_j) S_j^* f\rangle =\int
\psi\circ\sigma_j d\mu_{S_j^* f}, $$
with $S_j^*f =N^{-\delta/2} \chi_{D_j} f\circ \sigma_j$.

\smallskip

Iterating the relation \eqref{scaleFourier} one obtains
\begin{equation}\label{FmufSastar}
\hat\mu_f(t) = \sum_{a\in \cW_{k,A}} e^{it x(a)} \hat\mu_{S_a^* f}(\frac{t}{N^k}),
\end{equation}
where for $a=(a_1,\ldots,a_k)\in \cW_{k,A}$ we denote by $x(a)$ the
expression
\begin{equation}\label{xa}
x(a) = \frac{a_1}{N} + \frac{a_2}{N^2} + \cdots + \frac{a_k}{N^k}.
\end{equation}

\smallskip

As in \cite{Jor2}, we then obtain an approximation of the measure
$\mu_f$ with a family of combinations of Dirac measures in the 
following way.

\begin{cor}\label{approxDirac}
Let $\mu_f^{(k)}$ denote the measure
\begin{equation}\label{kDirac}
\mu_f^{(k)}(E) =\sum_{a\in \cW_{k,A}} \| S_a^* f \|^2 \delta_{a}(E),
\end{equation}
where $\delta_{a}$ is the Dirac measure supported at the rational
point $x(a)$ in $\Lambda_A$ whose terminating $N$-adic expansion is of the
form \eqref{xa}, for $$a=(a_1,\ldots,a_k)\in \cW_{k,A}.$$ 
The measures $\mu_f^{(k)}$ weakly converge to $\mu_f$, when considered
as functionals on the space of integrable functions $\psi$ on the real
line whose Fourier transform satisfies
\begin{equation}\label{tFconv}
 \int |t \hat\psi(t)| dt < \infty. 
\end{equation}
\end{cor}

\proof We show that, for all functions $\psi$ with \eqref{tFconv}, 
we have
$$ \lim_{k\to \infty} \int_{\Lambda_A} \psi \, d\mu_f^{(k)}
=\int_{\Lambda_A} \psi \,d\mu_f. $$
Passing to Fourier transforms, we have
$$ \int \psi \, d\mu_f^{(k)} - \int \psi \,d\mu_f = \int \hat \psi(t) (\hat
\mu_f^{(k)}(t)- \hat \mu_f(t)) \frac{dt}{2\pi}. $$
The Fourier transform of $\mu_f^{(k)}$ is clearly of the form
$$ \hat  \mu_f^{(k)}(t) =\sum_{a\in \cW_{k,A}} e^{it x(a)} \| S_a^*
f\|^2, $$
with $x(a)$ as in \eqref{xa}, and one can estimate as in \cite{Jor2}
$$ | \hat \mu_f(t) - \hat  \mu_f^{(k)}(t) |\leq |t| N^{-k}. $$
This gives
$$ \left| \int \psi \, d\mu_f^{(k)} - \int \psi \,d\mu_f \right| \leq 
\frac{N^{-k}}{2\pi} \int |t \hat\psi(t)| \, dt $$
which gives the weak convergence $\mu_f^{(k)}\to \mu_f$.
\endproof

\subsection{Sierpinski fractals}

There is another fractal object, besides the limit set $\Lambda_A$, 
that is naturally associated to an $N\times N$-matrix $A$ with entries 
in $\{ 0,1 \}$. This is a Sierpinksi fractal constructed in the following way.
Consider the square $\bS=[0,1]\times [0,1]$ and write points $(x,y)\in \bS$
in terms of the $N$-adic expansion 
$$ (x,y)=(\frac{x_1}{N}+\frac{x_2}{N^2}+\cdots+\frac{x_k}{N^k}+\cdots,
\frac{y_1}{N}+\frac{y_2}{N^2}+\cdots+\frac{y_k}{N^k}+\cdots), $$
with $(x_i,y_i)\in \{ 0,\ldots,N-1\}\times \{ 0,\ldots,N-1 \}=\fA^2$, for all $i\geq 1$. 
We then consider the subset $\bS_A\subset \bS$ given by
\begin{equation}\label{SAsieripinski}
\bS_A=\{ (x,y)\in S \,|\, A_{x_i,y_i}=1, \,\,  \forall i\geq 1\}.
\end{equation}
This is a Sierpinski fractal whose iterative construction starts by subdividing
the unit square $\bS$ into the $N^2$ subsquares of size $N^{-2}$ consisting of
points $(x,y)$ with first digits of the $N$-adic expansion equal to given
$(i,j)\in \fA^2$. One then keeps among these only those for which $A_{ij}=1$.
The procedure is then iterated by subdividing each of the remaining squares 
into $N^2$ subsquares of size $N^{-4}$ and keeping only those for which
the same condition $A_{ij}=1$ is satisfied, and so on.
At each step a square is of size $N^{-2k}$ is replaced by $D$ squares of
size $N^{-2(k+1)}$, where
\begin{equation}\label{Ddi}
D = \sum_{i=0}^{N-1} d_i, \ \ \ \text{ with } \ \ \ d_i =\# \{ j \,|\, A_{ij}=1 \}.
\end{equation}
These satisfy $d_i\leq N$ and $D\leq N^2$.
Thus, the Hausdorff dimension of the Sierpinski fractal $\bS_A$ is simply
\begin{equation}\label{dimHSA}
\dim_H (\bS_A)= \frac{\log D}{2\log N}. 
\end{equation}

One can then consider maps $\tau_{(i,j)}: \bS_A \to \bS_A$, for $(i,j)$ satisfying
$A_{ij}=1$, given by
\begin{equation}\label{sigmaijSA}
\tau_{(i,j)} (x,y)=(\tau_i(x),\tau_j(y))=(\frac{x+i}{N},\frac{x+j}{N}).
\end{equation}
Notice how, unlike the $\sigma_i$ acting on $\Lambda_A$ that we considered before, 
here the $\tau_{(i,j)}$ are {\em everywhere defined} on $\bS_A$.  Since we are only considering 
such maps for pairs $(i,j)$ with $A_{ij}=1$, it is clear that the image $(\tau_i(x),\tau_j(y))$ is 
still a point in $\bS_A$. The corresponding coding map $\tau:\bS_A \to \bS_A$ is given by
$$ \tau(x,y)=(\tau(x),\tau(y))=
(0.x_2\cdots x_k\cdots,0.y_2\cdots y_k\cdots) , $$
for $(x,y)=(0.x_1x_2\cdots x_k\cdots,0.y_1y_2\cdots y_k\cdots)$.

\begin{lem}\label{ODreplem}
The semibranching function system $\{ \tau_{(i,j)} \}$ for $(i,j)\in \fA^2$ with $A_{ij}=1$
determines a representation of the Cuntz algebra $O_D$ on the Hilbert space $L^2(\bS_A,\mu)$, 
with $\mu$ the Hausdorff measure of dimension $\delta=\dim_H(\bS_A)$ as in \eqref{dimHSA}.
\end{lem}

\proof Let $\Phi_{(i,j)}$ denote the Radon--Nikodym derivative of the measure $\mu$ with respect
to composition by $\tau_{(i,j)}$. Since $\tau_{(i,j)}$ is of the form \eqref{sigmaijSA}, we have
\begin{equation}\label{RDijmu}
\Phi_{(i,j)}(x,y) = \frac{d\mu\circ \tau_{(i,j)}}{d\mu} = N^{-2\delta}=\frac{1}{D}.
\end{equation}

We consider the operators $S_{(i,j)}$ and $S_{(i,j)}^*$ defined as in the general case
of a semibranching function system in the form
\begin{equation}\label{Sij}
S_{(i,j)} f = \chi_{R_{(i,j)}}   \, \cdot \, (\Phi_{(i,j)}\circ \tau)^{-1/2}  \, \cdot \,  f\circ \tau ,
\end{equation}
with $R_{i,j}\subset \bS_A$ the range of $\tau_{(i,j)}$.
The adjoint $S_{(i,j)}^*$ in the inner product of $L^2(\bS_A,\mu)$ is given by
$$ \langle S_{(i,j)} f, h \rangle =N^\delta  \int_{R_{(i,j)} } f\circ \tau \, h \, d\mu =
N^\delta \int_{\bS_A} f \, h\circ\tau_{(i,j)}\, \Phi_{ij} \, d\mu , $$
so that we get
\begin{equation}\label{Sijstar}
S_{(i,j)}^* h = \Phi_{ij}^{1/2} \, h\circ\tau_{(i,j)} = N^{-\delta} \, h\circ\tau_{(i,j)}. 
\end{equation}
Thus, one sees that
\begin{equation}\label{ODrels}
S_{(i,j)}^* S_{(i,j)} = 1,  \ \ \ \text{ and } \ \ \  \sum_{(i,j): A_{ij}=1} S_{(i,j)} S_{(i,j)}^* = 1,
\end{equation}
since $S_{(i,j)} S_{(i,j)}^*$ is the range projection given by multiplication by $\chi_{R_{(i,j)}}$.
Thus, the $S_{(i,j)}$ generate a representation of the Cuntz algebra $O_D$ on
$L^2(\bS_A,\mu)$.
\endproof

In particular, this means that one can apply to the Sierpinski set $\bS_A$
all the techniques for constructions of wavelets on fractals from representations 
of Cuntz algebras developed, for instance, in \cite{Br-Jo}, \cite{Br-Jo-Os},
\cite{DuJor}, \cite{Jor}, \cite{Jor2}, \cite{Jor3}, etc.

Notice then that we can embed the limit set $\Lambda_A$ inside the Sierpinski
fractal $\bS_A$ in the following way.

\begin{lem}\label{embedSA}
The map
\begin{equation}\label{LambdaASA}
\Xi: \Lambda_A \to \bS_A, \ \ \ \  \Xi(x) =(x,\sigma(x)).
\end{equation}
gives an embedding $\Lambda_A \hookrightarrow \bS_A$. 
\end{lem}

\proof A point $x=(x_1x_2\cdots x_n\cdots)$ in $\Lambda_A$ satisfies
$A_{x_i x_{i+1}}=1$. This means that the point
$$ (x,y)=(0.x_1x_2\cdots x_n \cdots, 0.x_2x_3\cdots x_{n+1}\cdots)=(x,\sigma(x)) $$
satisfies $A_{x_i,y_i}=A_{x_ix_{i+1}}=1$ for all $i\geq 1$, hence it is a point in $\bS_A$. 
The map $\Xi$ is clearly injective since it is the identity on the first coordinate. It is continuous
since the preimage of a clopen set $\bS_A(i_1\cdots i_k,j_i\cdots j_k)$ of $\bS_A$, given
by numbers with fixed first $k$ digits of the $N$-adic expansion, is either empty, or else,
when $j_r=i_{r+1}$ for $r=1,\ldots,k-1$, it is equal to the clopen set
$\Lambda_A(i_1,\ldots,i_k,j_k)$ of $\Lambda_A$.
\endproof

One can then use this embedding together with the representation of the
algebra $O_D$ on $L^2(\bS_A,\mu)$ to obtain an induced action of a 
Cuntz--Krieger algebra.

\begin{prop}\label{ODOAinduce}
The maps $\tau_{i,j}$ restricts to maps defined on domains $D_{i,j}\subset \Xi(\Lambda_A)$.
These determine a semibranching function system on $\Xi(\Lambda_A)$ which gives rise
to a representation of the algebra $O_{\tilde A}$, where the $D\times D$-matrix $\tilde A$
is given by
\begin{equation}\label{OtildeA}
\tilde A_{(i,j),(\ell,k)} = \delta_{j,\ell}\, A_{jk}.
\end{equation}
\end{prop}

\proof The condition that $\tau_{(i,j)}(x,\sigma(x))=(\tau_i(x),\tau_j(\sigma(x))$ is in $\Xi(\Lambda_A)$ 
determines the domain $D_{(i,j)}\subset \Xi(\Lambda_A)$ to be 
\begin{equation}\label{DijXi}
D_{(i,j)}= \{ (x,\sigma(x))\in \Xi(\Lambda_A)\,|\, \sigma_j\sigma(x)=\sigma\sigma_i(x) \} =\Xi(R_j).
\end{equation}
In fact, the condition that $A_{ij}=1$ implies that $R_j\subset D_i$ in $\Lambda_A$, so that
$\Xi(D_i \cap R_j)=\Xi(R_j)$. We identify the restriction of continuous functions on $\bS_A$ 
to $\Xi(\Lambda_A)$ with continuous functions on $\Lambda_A$ and we write equivalently, 
with a slight abuse of notation, $f(x,\sigma(x))$ or $f(x)$. One then sees that
$$ f(\tau_{(i,j)}(x,\sigma(x))) = f(\sigma_i(x)) \,\, \chi_{R_j}(x). $$
This induces an isometry on the Hilbert space $L^2(\Xi(\Lambda_A),\mu_s)$,
where $\mu_s$ is the Hausdorff measure of dimension $s=\dim_H(\Xi(\Lambda_A))$,
$$ \hat S_{(i,j)}^* f(x) = N^s \chi_{R_{ij}} (x) f(\sigma(x)), $$
since for a function $f(x,\sigma(x))$ on $\Xi(\Lambda_A)$ we have
$$ \chi_{R_{(i,j)}}(x,\sigma(x)) f(\sigma(x),\sigma^2(x)) = \chi_{R_{ij}} (x) f(\sigma(x)). $$
This has adjoint 
$$ \hat S_{(i,j)}^* f (x) = N^{-s}  \chi_{R_j}(x) f(\sigma_i(x)). $$
We then see obtain
$$ \hat S_{(i,j)} \hat S^*_{(i,j)} f (x) = \chi_{R_{ij}} (x)\,\, \chi_{R_j}(\sigma(x))\,\, f(\sigma_i \sigma(x))
 = \chi_{R_{ij}}(x) \,\, f(x) $$
so that we have the relation
$$ \sum_{(i,j)}  \hat S_{(i,j)} \hat S^*_{(i,j)} =1. $$
We also have
$$ \hat S^*_{(i,j)} \hat S_{(i,j)}  f (x) = \chi_{R_j}(x) \, \chi_{R_{ij}}(\sigma_i(x)) f(\sigma \sigma_i(x)) =
\chi_{R_j}(x)\,\, f(x). $$
Using the fact that
$$ \chi_{R_j}= \sum_k A_{jk} \chi_{R_{jk}}, $$
we then obtain the other relation in the form
$$ \hat S^*_{(i,j)} \hat S_{(i,j)} = \sum_k A_{jk}  \hat S_{(j,k)} \hat S^*_{(j,k)}. $$
These correspond to the Cuntz--Krieger relations for the matrix $\tilde A$ of \eqref{OtildeA}. 
\endproof

\section{Wavelets on fractals}\label{WaveFracSec}

A general construction of wavelets on self-similar fractals was 
described in \cite{Jon}, see also \cite{Bodin}. The cases considered 
there correspond, from the point of view of semibranching function
systems, to the case where the $\sigma_i$ are defined on all of $X$,
as in the case of the Cuntz algebra. To adapt these constructions of
wavelets to the main case we are interested in, which is the Cantor
sets $\Lambda_A$ introduced above, one can use the representation
of the Cuntz--Krieger algebra $O_A$ on $L^2(\Lambda_A,d\mu_A)$ that we 
considered in the previous sections, and again the Perron--Frobenius
theory for the non-negative matrix $A$.

\smallskip

We begin by recalling briefly how the construction of \cite{Jon}
works in the case of a semibranching function
system on a measure space $(X,\mu)$ where the $N$ maps
$\sigma_i$ are defined on all of $X$. In this case one considers the 
$(m+1)$-dimensional linear space $\fP^m$ of polynomials on $\R$ of
degree $\leq m$, and one denotes by $\fS_0$ the
linear subspace of $L^2(X,d\mu)$, generated by the
restrictions $P|_{\Lambda_A}$ of polynomials in $\fP^m$. 
Under the condition that $X$ preserves Markov's inequality 
(see \S 4 of \cite{Jon}), one knows that one still has $\dim
\fS_0=m+1$. One then considers the linear
subspace $\fS_1\subset L^2(X,d\mu)$ of functions
$f\in L^2(X,d\mu)$ that are $\mu$-almost everywhere on $R_i=\sigma_i(X)$
restrictions $P|_{R_i}$ of some polynomial $P\in \fP^m$. 
Clearly $\fS_0 \subset \fS_1$ and $\dim \fS_1 = N \dim \fS_0 = N
(m+1)$, and let $\phi^\ell$, for $\ell=1,\ldots,m+1$ be an orthonormal
basis for $\fS_0$. One then considers the orthogonal complement 
$\fS_1\ominus \fS_0$, with a fixed choice of an orthonormal basis 
$\psi^\rho$, for $\rho =1,\ldots,(N-1)(m+1)$. 
The functions $\phi^r$ and $\psi^\rho$ provide the mother wavelets.  
One then considers the family of linear subspaces $\fS_k$ of
$L^2(X,d\mu)$, of functions whose restriction to
each subset $\sigma_{i_1}\circ \cdots \circ \sigma_{i_k}(X)$, 
agrees $\mu$-almost everywhere with the restriction to the same set 
of a polynomial in $\fP^m$. These satisfy $\fS_0\subset \fS_1 \subset 
\cdots \fS_k\subset \cdots L^2(X,d\mu)$. Moreover, any function in
$L^2(X,d\mu)$ can be approximated by elements in 
$$ \fS_0 \oplus \bigoplus_{k\geq 0} (\fS_{k+1}\ominus \fS_k), $$
since in fact the polynomials of degree zero already suffice, as
they give combinations of characteristic functions of the sets
$\sigma_{i_1}\circ \cdots \circ \sigma_{i_k}(X)$. The wavelets are
then obtained in \cite{Jon} as
\begin{equation}\label{waveletsX}
\psi^\rho_a = \mu(\sigma_a(X))^{-1/2} \psi^\rho\circ \sigma_a^{-1}.
\end{equation}
for $a=(i_1,\ldots,i_k)$ and $\sigma_a=\sigma_{i_1}\circ \cdots \circ
\sigma_{i_k}$.

\smallskip

We show now how to adapt this construction to the case of the Cantor
sets $\Lambda_A$. For simplicity, we describe in full only the case
where one only considers locally constant functions, that is, where one
starts with the 1-dimensional space $\fP^0$. This is the case that is 
closest to the classical construction based on the Haar wavelets, \cite{dau}. 

\smallskip

On the space $\Lambda_A \subset [0,1]$, with the Hausdorff measure
$\mu=\mu_A$, let $\fS_k$ denote the linear
subspaces of $L^2(\Lambda_A,d\mu_A)$ obtained as above, starting from
the 1-dimensional space $\fP^0$. Let 
\begin{equation}\label{fellk}
\{ f^{\ell,k} \}_{k=0,\ldots,N-1; \ell=1,\ldots,d_k}, 
\end{equation}
with 
\begin{equation}\label{dkA}
d_k = \# \{ j\,|\, A_{kj}=1 \},
\end{equation}
be a family of locally constant functions on $\Lambda_A$ such 
that the support of $f^{\ell,k}$ is contained in $R_k$ and
\begin{equation}\label{intfellk}
\int_{R_k} \overline{f^{\ell,k}} f^{\ell',k} = \delta_{\ell,\ell'}. 
\end{equation}
We also require that
\begin{equation}\label{fellkint0}
\int_{R_k} f^{\ell,k} =0 , \ \  \forall \ell=1,\ldots,d_k.
\end{equation}

\begin{lem}\label{S2S1basis}
A family of functions $f^{\ell,k}$ as in \eqref{fellk}, satisfying
\eqref{intfellk} and \eqref{fellkint0}, can be constructed using
linear combinations of characteristic functions $\chi_{R_{kj}}$, where
$R_{kj}=\Lambda_{2,A}(kj)$.  The resulting $f^{\ell, k}$ give an
orthonormal basis of the space $\fS_2\ominus \fS_1$. 
\end{lem}

\proof To see that linear combinations of characteristic functions
$\chi_{R_{kj}}$ suffice to construct the functions $f^{\ell,k}$,
notice first that the $\chi_{R_{kj}}$ give an orthogonal basis 
for the space $\fS_2$, which is of dimension $\dim \fS_2 = \sum_k
d_k$. We then write the $f^{\ell,k}$ in the form
\begin{equation}\label{fellkcoeff}
f^{\ell,k} =\sum_j A_{kj} c^{\ell,k}_j \chi_{R_{kj}},
\end{equation}
where the conditions \eqref{intfellk} and \eqref{fellkint0} translate
into conditions on the coefficients of the form
\begin{equation}\label{intfellkcoeff}
\sum_j A_{kj} \bar c_j^{\ell,k} c_j^{\ell',k} p_{kj} = \delta_{\ell,\ell'},
\end{equation}
where we use the notation
\begin{equation}\label{pkj}
p_{kj} = \mu(R_{kj}) = N^{-2\delta_A} p_j,
\end{equation}
according to \eqref{muLambdaka}, where $p=(p_0,\ldots,p_{N-1})$ is
the Perron--Frobenius eigenvector $A p = r(A) p$ for the non-negative
matrix $A$. Similarly, the condition
\eqref{fellkint0} becomes
\begin{equation}\label{fellkint0coeff}
\sum_j A_{kj} c_j^{\ell,k} p_{kj}=N^{-2\delta_A} \sum_j A_{kj} c_j^{\ell,k} p_j =0,
\end{equation}
where we again use \eqref{pkj}.

Let us introduce the following notation for convenience.
Consider on $\C^{d_k}\subset \C^N$ the inner product
\begin{equation}\label{pinnprod}
\langle v,w \rangle_k := \sum_j A_{kj} \bar v_j w_j p_j.
\end{equation}
Let $\cV_k$ denote the orthogonal complement, in the inner product
\eqref{pinnprod} on $\C^{d_k}$ of the vector $u=(1,1,\ldots,1)$, 
and let $\{ c^{\ell,k}=(c^{\ell,k}_i) \}_{\ell=1,\ldots,d_k-1}$ be 
an orthonormal basis of $\cV_k$, in the inner product \eqref{pinnprod},
namely
\begin{equation}\label{ccoeff}
\langle c^{\ell,k},u\rangle_k =0, \ \ \ \text{ and } \ \ \ \langle c^{\ell,k}
c^{\ell',k}\rangle_k =\delta_{\ell,\ell'}. 
\end{equation}
Then for $c^{\ell,k}$ as above, one sees that the functions
\eqref{fellkcoeff} are an orthonormal family satisfying 
the conditions \eqref{intfellk} and \eqref{fellkint0}. 

The space spanned by the $f^{\ell,k}$ is contained in $\fS_2$
by construction. The condition \eqref{fellkint0} ensures that the 
functions $f^{\ell,k}$ are orthogonal to all the $\chi_{R_k}$, hence 
they are in $\fS_2\ominus \fS_1$. They span a space of dimension
$\sum_k (d_k -1)= \sum_k d_k -N = \dim \fS_2\ominus \fS_1$. 
\endproof

\begin{thm}\label{waveletsthm}
Suppose given an orthonormal basis $\{ f^{\ell,r} \}$ for $\fS_2
\ominus \fS_1$, constructed as in Lemma \ref{S2S1basis} above.
Consider then functions of the form
\begin{equation}\label{waveletspsi}
\psi^{\ell,r}_a =  S_a\, f^{\ell,r},
\end{equation}
for $a=(a_1,\ldots,a_k) \in \cW_{k,A}$, give an orthonormal basis
for the space $\fS_{k+1}\ominus \fS_k$
hence, for varying $a\in \cW_{k,A}$ and for all $k\geq 0$, 
they give an orthonormal basis of wavelets for $L^2(\Lambda_A,\mu)$.
\end{thm}

\proof We have shown in Lemma \ref{S2S1basis} that the functions
$f^{\ell,r}$, for $r=0,\ldots, N-1$ and $\ell =1, \ldots, d_r$, 
give an orthonormal basis of $\fS_2\ominus \fS_1$.
We then check that the functions $S_a f^{\ell,r}$ give an
orthonormal basis for $\fS_{k+1}\ominus \fS_k$. 
Since in the representation of $O_A$ on $L^2(\Lambda_A,d\mu_A)$
we have constant Radon--Nikodym derivatives 
$\Phi_{\sigma_i}=N^{-\delta_A}$, this gives
$$ S_j f = N^{\delta_A/2}\, \chi_{R_j}\,\, f\circ \sigma, $$ 
so that we then have
$$ S_a\, f^{\ell,r} = N^{\delta_A k/2}\, \chi_{\Lambda_{k,A}(a)} \,\,
f^{\ell,r} \circ \sigma^k . $$
For $a\in \cW_{k,A}$, we have
$$ \langle S_a f^{\ell,r}, S_{a'} f^{\ell',r'} \rangle =
N^{\delta_A k} \langle \chi_{R_a}\, f^{\ell,r}\circ \sigma^k,
\chi_{R_{a'}}\, f^{\ell',r'}\circ \sigma^k \rangle $$ $$ = N^{\delta_A k}
\delta_{a,a'} \int_{R_a} \overline{(f^{\ell,r}\circ \sigma^k)}\,
(f^{\ell',r'}\circ \sigma^k)\, d\mu, $$
where we write $R_a = \Lambda_{k,A}(a)$, for the range of $\sigma_a
=\sigma_{a_1}\circ \cdots \circ\sigma_{a_k}$. Notice then that we
have, for any function $f\in L^2(\Lambda_A,d\mu)$ and any $a\in
\cW_{k,A}$,
\begin{equation}\label{intRafsigma}
 \int_{R_a} f\circ \sigma^k \, d\mu = \int_{D_{a_k}} f \,
\frac{d\mu\circ \sigma_a}{d\mu} \, d\mu $$ $$ = N^{-\delta_A k}
\int_{D_{a_k}} f\, d\mu = N^{-\delta_A k} \sum_j A_{a_k j} \int_{R_j}
f\, d\mu. 
\end{equation}
Applied to the above this gives
$$ \langle S_a f^{\ell,r}, S_{a'} f^{\ell',r'} \rangle = 
\delta_{a,a'}\delta_{r,r'} A_{a_k r} \int_{R_r} \overline{f^{\ell,r}}
f^{\ell',r}\, d\mu =\delta_{a,a'}\delta_{r,r'}\delta_{\ell,\ell'}. $$
Thus the $S_a f^{\ell,k}$ form an orthonormal system.

The space spanned by these functions is contained in $\fS_{k+1}$ and
a counting of dimensions shows that it has the dimension of
$\fS_{k+1}\ominus \fS_k$. To see that the $S_a f^{\ell,k}$ are in fact
orthogonal to the elements of $\fS_k$ it suffices to compute
$$ \langle S_a f^{\ell,r}, \chi_{\Lambda_{k,A}(b)} \rangle =
\delta_{a,b} N^{\delta_A k}
\int_{R_a} f^{\ell,r}\circ \sigma^k\, d\mu $$ $$ = \delta_{a,b} \sum_j
A_{a_k,j} \int_{R_j} f^{\ell,r}\, d\mu = 
\delta_{a,b}A_{a_k,r}\int_{R_r} f^{\ell,r}\, d\mu = 0, $$
by \eqref{intRafsigma} and \eqref{fellkint0}. This shows that we
obtained an orthonornal basis of $\fS_{k+1}\ominus \fS_k$, hence a
wavelet system for $L^2(\Lambda_A,d\mu)$.
\endproof

It is useful to remark how the main difference in this case, as
opposed to the similar constructions given for instance in \cite{Jon}
that we mentioned above, is that here we need to start from an
orthonormal basis of $\fS_2\ominus \fS_1$ instead of $\fS_1\ominus
\fS_0$. This reflects the fact that our functions $\sigma_i$ are not
everywhere defined and, while the choice of an orthonormal basis for
$\fS_1\ominus \fS_0$ gives the needed information on the ranges $R_i$,
in order to control both the ranges and the domains $D_i$ one needs 
to go one step further before starting the induction that constructs 
the wavelets, and consider $\fS_2\ominus \fS_1$.
Thus, the wavelet decomposition of a 
function $f\in L^2(\Lambda_A, \mu)$ will be given by
\begin{equation}\label{wavedec}
f = \sum_{k=0}^{N-1}\sum_{\ell=1}^{d_k -1} \alpha_{\ell,k} \, f^{\ell,k} + \sum_{j=0}^\infty \,\, 
\sum_{a\in \cW_{j,A}} \sum_{(\ell,k)} \alpha_{\ell,k,a}\, S_a f^{\ell,k}.
\end{equation}

The more general case where one starts the wavelet construction 
from the linear space of polynomials $\fP^m$ with $m\geq 1$ can 
be done along the same lines as Lemma \ref{S2S1basis} and Theorem
\ref{waveletsthm}. We describe in the next section a different
approach to wavelets constructions based on the Ruelle transfer
operator for the coding map $\sigma$. This is closer to the point of
view developed in \cite{DuJor}.

\section{Ruelle transfer operator}

A more general version of the Perron--Frobenius operator associated to
the coding map $\sigma: \Lambda_A \to \Lambda_A$ is obtained by 
considering the Ruelle transfer operator. This depends on the choice 
of a potential function $W$, defined on $\Lambda_A$, and is defined as
\begin{equation}\label{Ruelle}
\cR_{\sigma,W}f(x) = \sum_{y\,:\,\sigma(y)=x} W(y)\, f(y).
\end{equation}

\begin{lem}\label{Ruelleadj}
If the function $W$ is real valued, one can describe
the operator $\cR_{\sigma,W}$ as the adjoint of the operator 
\begin{equation}\label{TWop}
T_W f (x)= N^{\delta_A} \, W(x)\, f(\sigma(x)).
\end{equation}
\end{lem}

\proof We have
$$ \langle T_W f ,h \rangle = \int_{\Lambda_A} N^{\delta_A} \, W(x)\,
\overline{f(\sigma(x))} h(x) \, d\mu(x) $$
$$ = \sum_i \int_{D_i} \overline{f(u)}\, W(\sigma_i(u)) h(\sigma_i(u))
d\mu(u), $$
using the fact that the Radon--Nikodym derivative
$d\mu\circ\sigma_i/d\mu=N^{-\delta_A}$. We then write the above as
$$ \sum_{i,j} A_{ij} \int_{R_j} 
\overline{f(u)}\, W(\sigma_i(u)) h(\sigma_i(u)) d\mu(u). $$ 
We also have
$$ \sum_{i,j} A_{ij} \chi_{R_j}(x) W(\sigma_i(x))
h(\sigma_i(x)) = \sum_i A_{ix_1} W(\sigma_i(x))h(\sigma_i(x)).  $$
Since the set of preimages of the point $x$ under the coding map is
given by
$$ \{ y\,|\, \sigma(y)=x \}=\bigcup_{i: A_{ix_1}=1} R_i, $$
we see that the above is in fact
$$ \sum_i A_{ix_1} W(\sigma_i(x))h(\sigma_i(x)) =
\sum_{y\,:\,\sigma(y)=x} W(y)\, f(y). $$
This shows that $\langle T_W f ,h \rangle =\langle f,
\cR_{\sigma,W}(h) \rangle$.
\endproof

We assume that the potential $W$ of the Ruelle transfer operator
satisfies the Keane condition, namely that it has non-negative real
values $W: \Lambda_A \to \R_+$, and satisfies
\begin{equation}\label{Keane}
\sum_{y:\sigma(y)=x} W(y) =1.
\end{equation}
Equivalently, this means
\begin{equation}\label{Keane2}
\sum_i A_{ix_1} W(\sigma_i(x)) =1.
\end{equation}

\subsection{Random processes}
In the same way as described in \cite{DuJor}, we relate here harmonic functions for the
Ruelle transfer operator, that is, functions satisfying $\cR_{\sigma,W} h=h$ to random
processes defined by transition probabilities for paths from a given point $x$ to the image
under the $\sigma_j$ and their iterates.

Let $A^t$ be the transpose of the matrix $A$. Then we have $a^t=(a_k,\ldots,a_1)\in \cW_{k,A}$ 
if and only if $a=(a_1,\ldots,a_k) \in \cW_{k,A^t}$.  We construct probability measures on the limit set 
$\Lambda_{A^t}$ that are related to fixed points of the Ruelle transfer operator for the coding
$\sigma: \Lambda_A \to \Lambda_A$. In the following we denote by $R_i$ and $D_i$, as before,
the ranges and domains of the maps $\sigma_i$ in $\Lambda_A$ and by $R_i^t$ and $D_i^t$ the
corresponding sets in $\Lambda_{A^t}$.

For a given potential $W$ on $\Lambda_A$ satisfying the Keane condition \eqref{Keane}, 
consider a function $x\mapsto P^W_x$,  for $x\in D_i\subset \Lambda_A$, where 
$P^W_x: \cB(\Lambda_{A^t}\cap R^t_i) \to \R_+$, is a
non-negative function on the Borel subsets of $\Lambda_{A^t}$ 
defined by assigning to the $\Lambda_{k,A^t}(a)$ the values
\begin{equation}\label{Pxvalues}
P^W_x(\Lambda_{k,A^t}(a)) = A_{a_1 x_1} W(\sigma_{a_1}(x))
W(\sigma_{a_2}\sigma_{a_1}(x))\cdots W(\sigma_{a_k}\cdots\sigma_{a_1}(x)),
\end{equation}
for $a^t=(a_k,\ldots,a_1)\in \cW_{k,A}$ and for $x\in D_{a_1} \subset \Lambda_A$.

\begin{lem}\label{measurePWx}
The assignment \eqref{Pxvalues}, for $x\in D_i \subset \Lambda_A$, defines a
measure on $R_i^t\subset \Lambda_{A^t}$.
\end{lem}

\proof Similarly, to the case of $\Lambda_A$ seen in \eqref{nucompatible}, to 
check that \eqref{Pxvalues} defines a measure one has to check the compatibility condition
\begin{equation}\label{Pxcompatible}
P^W_x (\Lambda_{k,A^t}(a)) = \sum_j A_{a_k j}^t\, P^W_x(\Lambda_{k+1,A}(a j)),
\end{equation}
for all $x\in D_{a_1}\subset \Lambda_A$. We have
$$ P^W_x(\Lambda_{k+1,A}(a j)) = A_{a_1, x_1} W(\sigma_{a_1}(x)) 
\cdots W(\sigma_{a_k}\cdots\sigma_{a_1}(x)) W(\sigma_j \sigma_{a_k}\cdots\sigma_{a_1}(x))) $$
Moreover, the Keane condition for $W$ on $\Lambda_A$ gives
$$ \sum_j A_{j a_k} W(\sigma_j \sigma_{a_k}\cdots\sigma_{a_1}(x))) =1, $$
so we obtain \eqref{Pxcompatible}.
\endproof

One can think of the values of the potential $W$ as defining a probability of transition, or walk, 
from $x$ to $\sigma_{a_1}(x)$, so that \eqref{Pxvalues} can be regarded as the probability of a
random walk from $x$ to $\sigma_{a_k}\cdots \sigma_{a_1}(x)$. 
We then see that the random process $P^W_x$ is related to the fixed points of the Ruelle transfer operator.

\begin{prop}\label{Pfixmeas}
The random process $x\mapsto P_x^W$ introduced above is related to fixed points
of the Ruelle transfer operator in the following ways.
\begin{enumerate}
\item Let $E\subset \Lambda_{A^t}$ be a shift invariant set $\sigma^{-1}(E)=E$. Then the 
function $x\mapsto P_x^W(E)$ is a fixed point of the Ruelle transfer operator with 
potential $W$ on $\Lambda_A$.
\item If the series
\begin{equation}\label{seriesPx}
h(x):= \sum_{k\geq 1} \sum_{a\in \cW_{k,A^t}} A_{a_1 x_1} W(\sigma_{a_1}(x))\cdots W(\sigma_{a_k}\cdots\sigma_{a_1}(x))
\end{equation}
converges, then the function $h(x)$ is a fixed point of the Ruelle transfer operator with 
potential $W$ on $\Lambda_A$.
\end{enumerate}
\end{prop}

\proof (1) We check that this condition is equivalent to the fixed point
condition under the Ruelle transfer operator. For a given set
$\Lambda_{k,A^t}(a)$, we have
$$ \cR_{\sigma,W}(P^W_x (\Lambda_{k,A^t}(a))) =\sum_{y:\sigma(y)=x} W(y)\,
P^W_y(\Lambda_{k,A^t}(a)) $$ $$ =
\sum_j A_{jx_1} W(\sigma_j(x))\, P_{\sigma_j(x)}^W(\Lambda_{k,A^t}(a)). $$
A shift invariant set $\sigma^{-1}(E)=E$ in $\Lambda_{A^t}$ satisfies
$$ \cup_{j,i: A_{ji}^t =1} \sigma_j(E\cap R_i) = E. $$
By construction of the measures $P_x^W$, we know that $P_x^W(\sigma_j(E\cap R^t_i))$ is 
non-trivial provides that $x\in D_j$, so that $A_{jx_1}=1$. Thus, for $\sigma^{-1}(E)=E$,
we have
$$ \cR_{\sigma,W}(P^W_x (E)) = \sum_j A_{ij} P^W_x(\sigma_j(E\cap R^t_i))=
P^W_x(\sigma^{-1}(E)) = P^W_x(E), $$
which shows that $P^W_x(E)$ is a fixed point for $\cR_{\sigma,W}$.

(2) Assuming that the series \eqref{seriesPx} converges, we have
$$ \cR_{\sigma,W} h (x) = \sum_{\sigma(y)=x} W(y) h(y)  =
\sum_j A_{jx_1} W(\sigma_j(x)) h(\sigma_j(x)) $$
$$ = \sum_j A_{jx_1}  W(\sigma_j(x)) \sum_k \sum_a A_{a_1 j} W(\sigma_{a_1}\sigma_j(x))  \cdots W(\sigma_{a_k}\cdots\sigma_{a_1}\sigma_j(x)) $$
$$ = \sum_k \sum_{b=ja\in W_{k+1,A^t}} A_{jx_1}   W(\sigma_j(x)) W(\sigma_{a_1}\sigma_j(x))  \cdots 
W(\sigma_{a_k}\cdots\sigma_{a_1}\sigma_j(x)). $$
This gives $\cR_{\sigma,W} h (x) =h(x)$.
\endproof

\subsection{A trigonometric example}

We give an example of a potential $W$ satisfying the Keane condition, constructed using trigonometric functions. 

\begin{lem}\label{Wtrigolem}
The function
\begin{equation}\label{Wtrigo}
W(x) = \frac{1}{N_1} \left( 1- \cos\left( \frac{2\pi N x}{N_1} \right) \right) ,
\end{equation}
with $N_1=\#\{ j\,:\, A_{jx_1}=1 \}$,
is a potential satisfying the Keane condition \eqref{Keane} on $\Lambda_A$.
\end{lem}

\proof First notice that we have
$$ \sum_{j=0}^{N-1} A_{j x_1} \exp\left( \frac{2\pi i N \sigma_j(x)}{N_1} \right) =0 , $$
since $\sigma_j(x)= (x+j)/N$ and the above becomes a sum over all the $N_1$-th roots of unity. 
It follows directly from this that the real valued trigonometric version also satisfies 
$$ \sum_{j=0}^{N-1} A_{j x_1} \cos\left( \frac{2\pi N \sigma_j(x)}{N_1} \right) =0, $$
from which it follows that the potential of \eqref{Wtrigo} satisfies
$$  \sum_{j=0}^{N-1} A_{j x_1} W(\sigma_j(x)) =1. $$
Moreover, the function $W(x)$ takes non-negative real values, so it gives a potential with 
the Keane condition.
\endproof

\section{Examples and applications}

\subsection{Hensley Cantor sets and continued fraction expansion}

In \cite{ManMar} the coding of geodesics on the modular curves
$X_\Gamma = \H/\Gamma$, for $\Gamma\subset {\rm PGL}_2(\Z)$ a finite index subgroup
and $\H$ the hyperbolic upper half plane, was related to a generalization of the
shift map of the continued fraction expansion  $T: [0,1] \times \P \to [0,1] \times \P$,
\begin{equation}\label{Tshift}
 T(x,s) = \left( \frac{1}{x} - \left[\frac{1}{x}\right] , \left( \begin{array}{cc}
 -[1/x] & 1 \\ 1 & 0 \end{array} \right) s \right),
\end{equation}
where $\P ={\rm PGL}_2(\Z)/\Gamma$ is the finite coset set. It was then shown in
\cite{Mar}, \cite{Mar2}, that the restriction of this dynamical system to the Hensley Cantor sets,
that is, those subsets $E_N \subset [0,1]$ of points that only contains digits 
$a_k \leq N$ in the continued fraction expansion, gives rise to a dynamical system
\begin{equation}\label{sigmaE}
 \sigma : E_N \times \P \to E_N \times \P ,
\end{equation} 
which can be identified with the coding map $\sigma: \Lambda_A \to \Lambda_A$ 
of a semibranching function system $\{ \sigma_i \}$ that determines a Cuntz--Krieger algebra
$O_A$. The case where $\Gamma = {\rm PGL}_2(\Z)$ recovers the Cuntz algebra $O_N$. 
 
In this setting, one considers the Ruelle transfer operator with potential (without Keane condition)
$$ W(x,s)= |T'(x,s)|^\beta $$
so that
$$ \cR_{T, W} f (x,s) = \sum_{T(y,t)=(x,s)} |T'(y,t)|^\beta f(y,t) $$ $$ = \sum_{n=1}^N \frac{1}{(x+n)^{2\beta}} f\left(\frac{1}{x+n},   \left(\begin{array}{cc} 0 & 1 \\ 1 & n \end{array}\right) s\right). $$
This can be written in the form
$$ \sum_{(n,t)} A_{(n,t),(x_1,s)} W(\sigma_{(n,t)}(x,s)) f(\sigma_{(n,t)}(x,s)), $$
where the matrix $A$ is defined by the condition
$$ A_{(n,t),(k,s)}= \left\{ \begin{array}{ll} 1 &  M_n s =t \\[2mm]
0 & \text{otherwise,} 
\end{array} \right. $$ 
where the matrix $M_n\in {\rm GL}_2(\Z)$, acting on the left on the coset $\P$, is
$$ M_n =\left(\begin{array}{cc} 0 & 1 \\ 1 & n
\end{array}\right). $$

The shift invariant measure $\mu_{N,\P}$ on $E_N\times \P$ constructed in \cite{Mar} using the
fixed point of the Ruelle transfer operator can then be also seen as in \cite{KeStaStr}
as KMS$_\beta$ state for the time evolution on the Cuntz--Krieger algebra $O_A$
given by
$$ \sigma_t (S_{(k,s)})=  W^{-it}  S_{(k,s)},  $$
where we identify $W^{-it}$, for fixed $t$, with an element in $C(E_N\times \P)$.
The KMS state is then of the form
$$ \varphi_\beta(S_a S_a^*) = \int_{E_N\times \P} f_a(x,s) \, d\mu_{N,\P}(x,s) , $$
for $a=((k_1,s_1),\ldots, (k_r,s_r))\in \cW_{r,A}$ and $f_a$ the element in
$C(E_N\times \P)$ that corresponds to $S_a S_a^*$.  
The Ruelle operator can correspondingly be written as
$$ \cR_{T,W} f = \sum_{(n,t)} S^*_{(n,t)} \, W f\,  S_{(n,t)} $$
in term of generators of the Cuntz--Krieger algebra. 

\subsection{Graph wavelets from Cuntz--Krieger algebras}

It was recently shown, see for instance \cite{Crov}, that the
crucial problem of spatial traffic analysis on networks can be
addressed using a form of wavelet analysis which is adapted 
to the topology of the network graph. These {\em graph wavelets}
are constructed as families  of functions $\Psi_{\alpha}(v)$ 
on the set of vertices  $V(G)$ of a given finite graph $G$, 
localized with respect to certain scaling indices $\alpha$, and
with the property that
\begin{equation}\label{intPsiG}
\int_{V(G)} \Psi_\alpha(v) d\mu(v) =0, \ \ \ \text{ and } \ \ \  
\int_{V(G)} \bar \Psi_\alpha(v) \Psi_{\alpha'}(v) d\mu(v) = \delta_{\alpha,\alpha'},
\end{equation}
where $\mu(v)$ is a given measure that weights the nodes of the network
with assigned probabilities. 
We show here how to construct families of graph wavelets using the representations of
Cuntz--Krieger algebras and the corresponding wavelets on $\Lambda_A$ constructed 
in \S \ref{WaveFracSec}. 

Let $G$ be a finite directed graph with no sinks. It is well known that one can associate
to such a graph a Cuntz--Krieger algebra in the following way. One considers 
a collection of projections $P_v$ associated to the vertices $v\in V(G)$ and a
collection of partial isometries $S_e$ associated to the oriented
edges $e\in E(G)$, with the relations
\begin{equation}\label{relGraphCK1}
P_v = \sum_{s(e)=v} S_e S_e^* 
\end{equation}
for all $v\in V(G)$, and
\begin{equation}\label{relGraphCK2} 
P_{r(e)} = S_e^* S_e,
\end{equation}
for all edges $e\in E(G)$. Assuming that the graph has no sinks, so that
all vertices are sources, one has $\sum_v P_v=1$ so that the isometries
$S_e$ satisfy the relation \eqref{CKrel2},
\begin{equation}\label{relGraphCK1a}
\sum_e S_e S_e^* =1 .
\end{equation}
Moreover, for $N=\# E(G)$, one defines the $N\times N$-matrix $A_{ee'}$ by
\begin{equation}\label{AGmatrix}
A_{ee'}= \left\{ \begin{array}{ll} 1 & r(e)=s(e') \\
0 & \text{otherwise.} \end{array}\right.
\end{equation}
Then the relation \eqref{relGraphCK2}  reads equivalently as
\begin{equation}\label{relGraphCK2a}
S_e^* S_e = \sum_{e': r(e)=s(e')} S_{e'} S_{e'}^* = \sum_{e'} A_{ee'} S_{e'} S_{e'}^*,
\end{equation}
which gives the other Cuntz--Krieger relation \eqref{CKrel1}.

As before, let $\Lambda_A$ be the limit set associated to the algebra $O_A$
of the graph $G$. Let $d_e =\#\{ e'\,|\, r(e)=s(e')\}= \#\{ e'\,|\, A_{ee'}=1\}$. Consider as in
\S \ref{WaveFracSec} the orthonormal family of functions $\{ f^{\ell, e} \}$ with
$e\in E(G)$ and $\ell =1, \ldots, d_e$. As we have seen in \S \ref{WaveFracSec} these
are the mother wavelets for the orthonormal basis of $L^2(\Lambda_A,\mu_A)$ given by
the functions $\{ S_a f^{\ell,e} \}$, for varying $a\in \cW_{k,A}$ and $k\in \N$. Here an
element $a=(e_1,\ldots,e_k)\in  \cW_{k,A}$ is a path in the graph $G$ of length $k$ starting at 
the vertex $s(e_1)$.  Here we use the same mother functions to construct a family
of graph wavelets. 

Recall from \S \ref{WaveFracSec} that the functions $f^{\ell,e}$ are constructed in
terms of a family $c^{\ell,e}=(c^{\ell,e}_{e'})$ of vectors satisfying
\begin{equation}\label{celle1}
\sum_{e'} A_{e,e'} \bar c^{\ell,e}_{e'} c^{\ell',e}_{e'} p_{ee'} = \delta_{e,e'} , 
\end{equation}
where $p_{ee'}=\mu(R_{ee'})=N^{-2\delta_A} p_{e'}$ and
\begin{equation}\label{celle2}
\sum_{e'} A_{e,e'} c^{\ell,e}_{e'} p_{e'} = 0 . 
\end{equation}
Upon rescaling the coefficients $c^{\ell,e}_{e'}$ by a factor $N^{\delta_A}$,
we obtain a family satisfying \eqref{celle2} and with \eqref{celle1} replaced
by the similar
\begin{equation}\label{celle}
\sum_{e'} A_{e,e'} \bar c^{\ell,e}_{e'} c^{\ell',e}_{e'} p_{e'} = \delta_{e,e'} ,
\end{equation}
where we keep the same notation for these rescaled coefficients. The $p_e$ are
the components of the Perron--Frobenius eigenvector $A p= r(A) p$.

After fixing a choice of a base vertex $v_0 \in E(G)$, we define 
a measure on the set of vertices of the graph by $\mu_{G,v_0}(v_0)=0$ and
\begin{equation}\label{muGv}
\mu_{G,v_0}(v) := p_{e_1} \cdots p_{e_k},
\end{equation}
where $e_1 \cdots e_k$ is the shortest path in the graph $G$ starting at $v_0$ and 
ending at $v$. This means that we are considering a random walk on the graph 
starting at $v_0$, where at the first step one has probability $p_e$ of moving to the
nearby vertex $r(e)$ and probability zero of remaining at $v_0$. The measure
\eqref{muGv} gives the probability of reaching at time $k$ one of the vertices that
are $k$ steps away from $v_0$.

In addition to fixing the base vertex $v_0$, we also fix a choice of an edge $e_0$
with $r(e_0)=v_0$. We then define functions 
\begin{equation}\label{waveG1}
\Psi_{\ell}(v) =\left\{ \begin{array}{ll} c^{\ell,e_0}_{e'} & v=r(e'),\, \, v_0=r(e_0)=s(e') \\[2mm] 
0 & \text{otherwise.}
\end{array}\right.
\end{equation}
These satisfy
\begin{equation}\label{waveGint1}
\int_{V(G)} \Psi_\ell(v) d\mu_{G,v_0}(v) =
\sum_{e'} A_{e_0 e'} c^{\ell,e_0}_{e'} p_{e'} = 0 
\end{equation}
and
\begin{equation}\label{waveGint2}
\int_{V(G)} \bar\Psi_\ell (v) \Psi_{\ell'}(v)d\mu_{G,v_0}(v) =
\sum_{e'} A_{e_0 e'} \bar c^{\ell,e_0}_{e'} c^{\ell',e_0}_{e'}  p_{e'} = \delta_{\ell,\ell'}.
\end{equation}
We then extend this to a family $\Psi_{\ell_1,\ldots,\ell_k}(v)$, where we
consider paths $a=(e_1,\ldots, e_k)\in \cW_{k,A}$ of length $k$ in the graph
starting at $v_0$, with $\ell_i =1, \ldots, d_{e_i}$. We set
\begin{equation}\label{waveG2}
\Psi_{\ell_1,\ldots,\ell_k}(v) =\left\{ \begin{array}{ll} c^{\ell_1,e_0}_{e_1}  c^{\ell_2,e_1}_{e_2}\cdots
c^{\ell_k,e_{k-1}}_{e_k} & v=r(e_k),\,\, v_0=s(e_1) \\[2mm]
0 & \text{otherwise.}
\end{array}\right.
\end{equation}
These again satisfy
\begin{equation}\label{waveGint1a}\begin{array}{c}
\displaystyle{\int_{V(G)} \Psi_{\ell_1,\ldots,\ell_k}(v) d\mu_{G,v_0}(v)} =\\[3mm]
\displaystyle{\sum_{(e_1,\ldots,e_k)} A_{e_0 e_1}\cdots A_{e_{k-1} e_k} 
c^{\ell_1,e_0}_{e_1}  c^{\ell_2,e_1}_{e_2}\cdots
c^{\ell_k,e_{k-1}}_{e_k}  p_{e_1} \cdots p_{e_k} =0}. \end{array}
\end{equation}
This vanishes since already $\sum_{e_k} A_{e_{k-1} e_k} c^{\ell_k,e_{k-1}}_{e_k}  p_{e_k}=0$. 
Moreover, they satisfy
\begin{equation}\label{waveGint2a}\begin{array}{c}
\displaystyle{\int_{V(G)} \bar\Psi_{\ell_1,\ldots,\ell_k} (v) \Psi_{\ell'_1,\ldots,\ell'_k}(v)d\mu_{G,v_0}(v) }=
\\[3mm] \displaystyle{
\sum_{(e_1,\ldots,e_k)} A_{e_0 e_1} \cdots A_{e_{k-1} e_k}
\bar c^{\ell_1,e_0}_{e_1}  \cdots \bar c^{\ell_k,e_{k-1}}_{e_k} 
c^{\ell_1',e_0}_{e_1} \cdots c^{\ell_k',e_{k-1}}_{e_k}  p_{e_1}\cdots p_{e_k} }=\\[3mm]
\delta_{\ell_1,\ell_1'}\cdots \delta_{\ell_k,\ell_k'}. \end{array}
\end{equation}
The functions $\Psi_{\ell_1,\ldots,\ell_k}$, for $k\geq 1$, constructed in this way, 
are supported on concentric regions $\cU_k(v_0)$ made of vertices at a distance $k$ from a
chosen base vertex $v_0$. Unlike other types of graph wavelets constructions where the 
functions are constant on such concentric regions $\cU_k(v_0)$ and average to zero 
over different  $k$, the ones we obtain here are supported on a single $\cU_k(v_0)$ 
with zero average. In terms of traffic analysis on networks, while one type of graph
wavelets may be more suitable in analyzing radial propagation from a vertex, the other
may be preferable for directional propagation away from a chosen vertex.

In \cite{ConsMa2}, \cite{CMRV} one considered, in the setting of Mumford curves 
with p-adic Schottky uniformization, the Cuntz--Krieger algebras associated to 
the finite graphs with no sinks obtained from the action of a p-adic Schottky 
group on the subtree of the Bruhat--Tits tree spanned by geodesics with boundary 
points on the limit set in $\P^1(\Q_p)$. In that context it would be interesting to
compare the wavelet constructions described in this paper with the p-adic wavelet
theory (see for instance \cite{AlKo}).

\medskip

{\bf Acknowledgments.}
We thank the Max-Planck Institut f\"ur Mathematik 
in Bonn, where part of this work was done, for the hospitality and
support and for the excellent working conditions. The first author is 
partially supported by NSF grant DMS-0651925.


\end{document}